\documentclass[11pt]{article}
\usepackage{amsfonts}
\usepackage{euscript,amstext,a4}
\usepackage{mathrsfs}
\usepackage{indentfirst}
\usepackage{enumerate}
\usepackage{amssymb}
\usepackage{amsmath}
\usepackage{hyperref}
\usepackage{amscd}
\numberwithin{equation}{section}
\newtheorem{thm}{Theorem}[section]
\numberwithin{thm}{section}
\newtheorem{defn}[thm]{Definition}
\newtheorem{prop}[thm]{Proposition}
\newtheorem{rem}[thm]{Remark}
\newtheorem{lem}[thm]{Lemma}

\DeclareMathOperator{\tr}{Tr}
\DeclareMathOperator{\SF}{sf}
\DeclareMathOperator{\APS}{APS-ind}
\DeclareMathOperator{\topo}{t-ind}

\def\R{{\mathbb{R}}}

\title{Rigidity and Vanishing Theorems \\ on ${\mathbb{Z}}/k$ Spin$^c$ manifolds}
\author{Bo LIU\footnote{Chern Institute of Mathematics \& LPMC, Nankai University, Tianjin
300071, P. R. China. (boliumath@mail.nankai.edu.cn)} \ \ and\  \
Jianqing YU\footnote{Chern Institute of Mathematics \& LPMC, Nankai
University, Tianjin 300071, P. R. China.
(jianqingyu@gmail.com)}}
\date{}

\begin{document}

\maketitle

\begin{abstract}
In this paper, we first establish an $S^1$-equivariant index theorem for Spin$^c$ Dirac operators on $\mathbb{Z}/k$  manifolds,
then combining with the methods developed by Taubes \cite{MR998662} and Liu-Ma-Zhang \cite{MR1870666,MR2016198}, we extend  Witten's
rigidity theorem to the case of $\mathbb{Z}/k$ Spin$^c$ manifolds. Among others, our results resolve a conjecture
of Devoto \cite{MR1405063}.
\end{abstract}

\section{Introduction}

In \cite{MR970288}, Witten derived a series of elliptic operators on the free loop space $\mathcal{L}M$ of
a spin manifold $M$. In particular, the index of the formal signature operator on loop space turns out to be exactly
the elliptic genus constructed by Landweber-Stong \cite{MR948178} and Ochanine \cite{MR895567} in a topological way. Motivated by
physics, Witten conjectured that these elliptic operators should be rigid with respect to the circle action.

This conjecture was first proved by Taubes \cite{MR998662} and Bott-Taubes \cite{MR954493}.
See also \cite{MR981372} and \cite{MR1048541} for other interesting cases.
By the modular invariance property, Liu (\cite{MR1331972,MR1396769}) presented a simple
and unified proof of the above conjecture as well as various further generalizations.
In particular, several new vanishing theorems were established in \cite{MR1331972,MR1396769}.
Furthermore, on the equivariant Chern character level, Liu and Ma (\cite{MR1756105,MR1969037})
generalized Witten's rigidity theorem
to the family case, and also obtained several vanishing theorems for elliptic genera.
In \cite{MR1870666,MR2016198}, inspired by \cite{MR998662}, Liu, Ma and Zhang established the corresponding family
rigidity and vanishing theorems on the equivariant $K$-theory level.

In \cite{MR1993996}, Zhang established an equivariant index theorem
for circle actions on $\mathbb{Z}/k$ spin manifolds and pointed out
that by combining with the analytic arguments developed in
\cite{MR2016198}, one can prove an extension of Witten's rigidity
theorem to $\mathbb{Z}/k$ spin manifolds. The purpose of this paper
is to extend the result of \cite{MR1993996} to $\mathbb{Z}/k$
Spin$^c$ manifolds and then establish Witten's rigidity theorem for
$\mathbb{Z}/k$ Spin$^c$ manifolds. Recall that a
$\mathbb{Z}/k$ manifold $X$ is a smooth manifold with boundary
$\partial X$ which consists of $k$ disjoint pieces, each of which is
diffeomorphic to a given closed manifold $Y$ (cf. \cite{MR0350748}). It is
interesting that for a  Dirac operator $D$ on a
$\mathbb{Z}/k$ manifold, the $\APS(D)\ {\rm mod}\ k\mathbb{Z}$ determines a
topological invariant in $\mathbb{Z}/k\mathbb{Z}$, where $\APS(D)$
is the index of $D$ which is imposed the boundary condition of
Atiyah-Patodi-Singer type \cite{MR0397797}. Freed and Melrose
\cite{MR1144425} proved a {\rm mod} $k$ index theorem,
\begin{equation}\label{fm1}
\APS(D)\ {\rm mod}\ k\mathbb{Z}=\topo(D)\,,
\end{equation}
giving the $\APS(D)\ {\rm mod}\ k\mathbb{Z}$ a purely topological interpretation.

Assume that $X$ is a $\mathbb{Z}/k$ manifold which admits
a $\mathbb{Z}/k$ circle action (cf. Section \ref{sec2.2}). Let $D$ be a Dirac operator on $X$ which commutes with the
circle action. Let $R(S^1)$ denote the representation ring of $S^1$.
The equivariant topological index of $D$ is defined by Freed and Melrose \cite{MR1144425}
as an element of $\mathbb{Z}/k\mathbb{Z}\otimes R(S^1)$,
and we denote it by $\topo_{S^1}(D)$.
Then there exist $R_n\in \mathbb{Z}/k\mathbb{Z}$ such that
\begin{equation}
\topo_{S^1}(D)=\sum_{n\in \mathbb{Z}}R_n\otimes[n]\,,
\end{equation}
where by $[n]$ ($n\in\mathbb{Z}$) we mean the one
dimensional complex vector space on which $S^1$ acts as
multiplication by $g^n$ for a generator $g\in S^1$.

On the other hand, by applying the equivariant index theorem for $\mathbb{Z}/k$ manifolds
established by Freed and Melrose in \cite{MR1144425},
one gets for $n\in \mathbb{Z}$,
\begin{equation}\label{fm2}
R_n=\APS(D,n)\ {\rm mod}\ k\mathbb{Z}\,.
\end{equation}
See \eqref{ly5} for the definition of $\APS(D,n)$.

The Dirac operator $D$ on $X$ is said to be {\em rigid in $\mathbb{Z}/k$ category}
for the circle action if its equivariant topological index
$\topo_{S^1}(D)$ verifies that for $n\in\mathbb{Z}$, $n\neq 0$, one has
\begin{equation}\label{uy20}
R_n=0\ \  \text{in}\ \mathbb{Z}/k\mathbb{Z}.
\end{equation}
Furthermore, we say $D$ has {\em vanishing property in $\mathbb{Z}/k$ category} if
its equivariant topological index $\topo_{S^1}(D)$ is identically zero, i.e., \eqref{uy20}
holds for any $n\in\mathbb{Z}$.

In \cite{MR1405063}, Devoto introduced what he called mod $k$ elliptic genus for $\mathbb{Z}/k$
spin manifolds as an $S^1$-equivariant topological index in the sense of \cite{MR1144425} of some
twisted Dirac operator and conjectured that this mod $k$ elliptic genus is {\em rigid in $\mathbb{Z}/k$ category}.
In this paper, following the suggestion in \cite[Remark 1]{MR1993996},
we present a proof of Devoto's conjecture. Moreover, we
establish our results for $\mathbb{Z}/k$ Spin$^c$ manifolds,
thus generalizing \cite[Theorems A and B]{MR1396769} to the case of
$\mathbb{Z}/k$ Spin$^c$ manifolds.

Our proof of these rigidity results consists of two steps. In step 1 (Sections \ref{sec2} and \ref{sec3}), we extend the $\mathbb{Z}/k$ equivariant index theorem of Zhang \cite{MR1993996} to the Spin$^c$ case.
In step 2 (Sections \ref{sec4} and \ref{sec5}), using the {\rm mod} $k$ localization index theorem established in step 1 and modifying the process in \cite{MR1870666,MR2016198}, we prove the main results of this paper.

This paper is organized as follows. In Section \ref{sec2}, we state an $S^1$-equivariant index theorem for Spin$^c$ Dirac operators on $\mathbb{Z}/k$ manifolds (cf. Theorem \ref{liu28}). As an application, we extend Hattori's vanishing theorem \cite{MR508087} to the case of $\mathbb{Z}/k$ almost complex manifolds. In Section \ref{sec3}, we
prove the $S^1$-equivariant index theorem stated in Section \ref{sec2}. In Section \ref{sec4}, we prove our main results (cf. Theorem \ref{maintheorem}), the rigidity and vanishing theorems for $\mathbb{Z}/k$ Spin$^c$ manifolds, which generalize \cite[Theorems A and B]{MR1396769}.
When applied to $\mathbb{Z}/k$ spin manifolds, our results
resolve a conjecture of Devoto \cite{MR1405063}.
Section \ref{sec5} is devoted to a proof of the recursive
formula which has been used in Section \ref{sec4} in the proof of our main results.

\section{Spin$^c$ Dirac operators and a mod $k$ localization formula}\label{sec2}

In this section, for a ${\mathbb{Z}}/k$ manifold which admits
a nontrivial ${\mathbb{Z}}/k$ circle action, we state a mod $k$ localization formula
for $S^1$-equivariant Spin$^c$ Dirac operators, whose proof will be
given in Section \ref{sec3}.
As an application, we deduce the rigidity and vanishing property for several Dirac operators
on a ${\mathbb{Z}}/k$ almost complex manifold.
In particular, we extend Hattori's vanishing theorem \cite{MR508087}
to the case of $\mathbb{Z}/k$ almost complex manifolds.

This section is organized as follows. In Section \ref{sec2.1}, we review the construction of Spin$^c$ Dirac operators on $\mathbb{Z}/k$ manifolds and the Atiyah-Patodi-Singer boundary problems. In Section \ref{sec2.2}, we recall
the circle actions on $\mathbb{Z}/k$ manifolds and present a variation formula for the indices
of these boundary problems. In Section \ref{sec2.3}, we state the
mod $k$ localization formula for ${\mathbb{Z}}/k$ circle actions. As an application, in Section \ref{sec2.4}, we extend Hattori's vanishing theorem \cite{MR508087} to the case of $\mathbb{Z}/k$ almost complex manifolds.

\subsection{Spin$^c$ Dirac operators on ${\mathbb{Z}}/k$ manifolds}\label{sec2.1}

We first recall the definition of $\mathbb{Z}/k$ manifolds introduced by Morgan and Sullivan (cf. \cite{MR0350748}).

\begin{defn}{\rm (cf. \cite[Definition 1.1]{MR1993996})}\label{def1}
A compact $\mathbb{Z}/k$ manifold is a compact manifold
$X$ with boundary $\partial X$, which admits a decomposition
$\partial X=\cup_{i=1}^k(\partial X)_i$ into $k$ disjoint manifolds
and $k$ diffeomorphisms $\pi_i:(\partial X)_i\rightarrow Y$ to a
closed manifold $Y$.
\end{defn}

Let $\pi:\partial X\rightarrow Y$ be the induced map. In what
follows, as in \cite{MR1993996}, we will call an object $\alpha$ (e.g., metrics, connections,
etc.) of $X$ a $\mathbb{Z}/k$-object if there will be a
corresponding object $\beta$ on $Y$ such that $\alpha|_{\partial
X}=\pi^*\beta$.

We point out here that in this paper when consider the
topological objects (e.g., cohomology, characteristic classes, $K$
group, etc.) on a $\mathbb{Z}/k$ manifold $X$, we always regard $X$
as a quotient space obtained by identifying each of the $k$ disjoint
pieces of the boundary $\partial X$. Then $X$ has the homotopy type
of a CW complex, which implies that the first Chern class $c_1$ induces a
1-to-1 correspondence between the equivalence classes of the complex
line bundles over $X$ and the elements of $H^2(X;\mathbb{Z})$.
As will be seen, this is essential in our proof.

We make the assumption that $X$ is $\mathbb{Z}/k$
oriented and of dimension $2l$.

Let $V$ be a $\mathbb{Z}/k$ real vector bundle over $X$ which is of dimension $2p$
and $\mathbb{Z}/k$ oriented.
Let $L$ be a $\mathbb{Z}/k$ complex line bundle over $X$ with the property that the
vector bundle $U=TX\oplus V$ satisfies $\omega_2(U)=c_1(L)\mod (2)$,
 where $\omega_2$ denotes the second Stiefel-Whitney
class, and $c_1$ denotes the first Chern class. Then the $\mathbb{Z}/k$ vector
bundle $U$ has a $\mathbb{Z}/k$ ${\rm Spin}^c$-structure.

Let $g^{TX}$ be a $\mathbb{Z}/k$ Riemannian metric on $X$. Let $g^{T\partial
X}$ be its restriction on $T\partial X$.  Let $\epsilon_0>0$ be less than the injectivity radius of $g^{TX}$.
We use the inward geodesic flow to identify a neighborhood of the boundary with the collar
$[0,\epsilon_0)\times\partial X$.
 We assume that $g^{TX}$ is
of product structure near $\partial X$. That is, there is an open
neighborhood $\mathscr{U}_{\epsilon}=[0,\epsilon)\times\partial X$ of
$\partial X$ in $X$ with $0<\epsilon\leq\epsilon_0$ such that one has
the orthogonal splitting on ${\mathscr U}_{\epsilon}$,
\begin{equation}\label{liu1}
g^{TX}|_{{\mathscr U}_{\epsilon}}=dr^2 \oplus \pi^*_{\epsilon}g^{T\partial
X},
\end{equation}
where $\pi_{\epsilon}:[0,\epsilon)\times\partial X\rightarrow \partial
X$ is the obvious projection onto the second factor.

Let $\nabla^{TX}$ be the Levi-Civita connection on $(TX, g^{TX})$. Then $\nabla^{TX}$ is a $\mathbb{Z}/k$ connection.

Let $W$ be a $\mathbb{Z}/k$ complex vector bundle over $X$ with a $\mathbb{Z}/k$ Hermitian metric
$g^W$. Let $\nabla^W$ be a $\mathbb{Z}/k$ Hermitian connection on $W$ with respect to $g^W$.
We make the assumption that $g^W$ and $\nabla^W$ are both of product structure near $\partial X$.
That is, over the open neighborhood $\mathscr{U}_{\epsilon}$ of
$\partial X$, one has
\begin{equation}\label{yy1}
g^W|_{{\mathscr U}_{\epsilon}}=\pi^*_{\epsilon}(g^W|_{\partial X}),\quad\text{and}
\quad\nabla^W|_{{\mathscr U}_{\epsilon}}=\pi^*_{\epsilon}(\nabla^W|_{\partial X}).
\end{equation}

Let $g^V$ (resp. $g^L$) be a $\mathbb{Z}/k$ Euclidean (resp. Hermitian) metric on $V$ (resp. $L$), and
$\nabla^V$ (resp. $\nabla^L$) be a corresponding $\mathbb{Z}/k$ Euclidean (resp. Hermitian)
connection on $V$ (resp. $L$). We make the assumption that
$g^V$, $\nabla^V$, $g^L$, $\nabla^L$ are of product structure near $\partial X$ (cf. (\ref{yy1})).

By taking $\epsilon>0$ sufficiently small, one can always find the
metrics $g^{TX}$, $g^W$, $g^V$, $g^L$ and the connections $\nabla^W$, $\nabla^V$, $\nabla^L$
verifying the above assumptions.

The Clifford algebra bundle $C(TX)$ is the bundle of Clifford algebras over $X$ whose fibre at $x\in X$ is the
Clifford algebra $C(T_xX)$ (cf. \cite{MR1031992}). Let $C(V)$ be the Clifford algebra bundle of
$(V,g^V)$.

Let $S(U,L)$ be the fundamental complex spinor bundle for
$(U,L)$ (cf. \cite[Appendix D]{MR1031992}). We denote by $c(\cdot)$ the
Clifford action of $C(TX)$, $C(V)$ on $S(U,L)$.
Let $\{e_i\}_{i=1}^{2l}$ (resp. $\{f_j\}_{j=1}^{2p}$) be an
oriented orthonormal basis of $(TX, g^{TX})$ (resp. $(V, g^V)$).
There are two canonical ways to consider $S(U,L)$ as a
$\mathbb{Z}_2$-graded vector bundle. Let
\begin{equation}\label{liu3}
\begin{split}
\tau_s&=(\sqrt{-1})^l c(e_1)\cdots c(e_{2l}),
\\
\tau_e&=(\sqrt{-1})^{l+p}c(e_1)\cdots c(e_{2l})c(f_1)\cdots
c(f_{2p})
\end{split}
\end{equation}
be two involutions of $S(U,L)$. Then $\tau_s^2=\tau_e^2=1$. We
decompose $S(U,L)=S_+(U,L)\oplus S_-(U,L)$ corresponding to $\tau_s$
(resp. $\tau_e$) such that $\tau_s|_{S_{\pm}(U,L)}=\pm1$ (resp.
$\tau_e|_{S_{\pm}(U,L)}=\pm1$).

In the remaining part of this paper, we always fix an involution $\tau$ on $S(U,L)$, either $\tau_s$ or $\tau_e$,
without further notice.

Let $\nabla^{S(U,L)}$ be the
Hermitian connection on $S(U,L)$ induced by
$\nabla^{TX}\oplus\nabla^V$ and $\nabla^L$ (cf. \cite[Appendix D]{MR1031992}). Then $\nabla^{S(U,L)}$
preserves the $\mathbb{Z}_2$-grading of $S(U,L)$. Let $\nabla^{S(U,L)\otimes W}$ be the Hermitian connection on
$S(U,L)\otimes W$ obtained from the
tensor product of $\nabla^{S(U,L)}$ and $\nabla^W$.

\begin{defn}\label{diracdef}
The twisted Spin$^c$ Dirac operator $D^{X}$ on $S(U,L)\otimes W$ over $X$ is defined by
\begin{equation}\label{d1}
D^X=\sum_{i=1}^{2l}c(e_i)\nabla_{e_i}^{S(U,L)\otimes W}:
\Gamma (X,S(U,L)\otimes W) \longrightarrow \Gamma (X,S(U,L)\otimes W).
\end{equation}
Denote by $D^X_{\pm}$ the restrictions of $D^X$ on
$\Gamma (X,S_\pm(U,L)\otimes W)$.
\end{defn}

By \cite{MR1031992}, $D^X$ is a formally self-adjoint operator. To get an elliptic operator, we
impose the boundary condition of Atiyah-Patodi-Singer type \cite{MR0397797}.

We first recall the canonical boundary operators (cf. \cite[(1.4)]{MR1870658}).
For a first order differential operator $D:\Gamma(S(U,L)\otimes W)\longrightarrow\Gamma(S(U,L)\otimes W)$ on $X$,
if there exists $\epsilon>0$ sufficient small such that the following identity holds on $\mathscr{U}_{\epsilon}$,
\begin{equation}\label{liu8}
D=c\left(\frac{\partial}{\partial r}\right)\left(\,\frac{\partial}{\partial r}+B\right)\,,
\end{equation}
with $B$ independent of $r$, then we
will call $B$ the canonical boundary operator associated to $D$.
When there is no confusion, we will also use $B$ to denote its
restriction on $\Gamma (X,S(U,L)\otimes W)|_{\partial X}$.

We then recall the Atiyah-Patodi-Singer projection associated to a boundary operator
(cf. \cite{MR0397797}). Assume temporarily that $B:\Gamma (X,S(U,L)\otimes W)|_{\partial X} \longrightarrow
\Gamma (X,S(U,L)\otimes W)|_{\partial X}$ is a first order formally self-adjoint elliptic differential
operator  on $\partial X$. For any $\lambda\in {\rm Spec}\,(B)$, the spectrum of $B$,
let $E_\lambda$ be the eigenspace corresponding to $\lambda$. For $a\in \mathbb{R}$, let $P_{\geq a}$
be the orthogonal projection from the $L^2$-completion of $\Gamma (X,S(U,L)\otimes W)|_{\partial X}$ onto
$\oplus_{\lambda\geq a}E_{\lambda}$. We call the particular projection
$P_{\geq 0}$ the Atiyah-Patodi-Singer projection associated to $B$ to emphasize its role in \cite{MR0397797}.
If we assume in addition that $B$ preserves the
$\mathbb{Z}_2$-grading of $\Gamma (X,S(U,L)\otimes W)|_{\partial X}$, and let $B_{\pm}$ be the restrictions of
$B$ on $\Gamma (X,S_{\pm}(U,L)\otimes W)|_{\partial X}$, then we will restrict $P_{\geq a}$
on the $L^2$-completions of $\Gamma (X,S_{\pm}(U,L)\otimes W)|_{\partial X}$ and denote them by $P_{\geq a,\pm}$.

Let $e_{1}=\frac{\partial}{\partial r}$ be the inward unit normal vector field perpendicular to $\partial X$.
Let $e_2,\cdots,e_{2l}$ be an oriented orthonormal basis of $T\partial X$ so that
$e_1,e_2,\cdots,e_{2l}$ is an oriented orthonormal basis of $TX|_{\partial X}$. Then
using parallel transport with respect to $\nabla^{TX}$ along the unit speed geodesics perpendicular
to $\partial X$, $e_1,e_2,\cdots,e_{2l}$ forms an oriented orthonormal basis of $TX$ over $\mathscr{U}_{\epsilon}$.

\begin{defn}\label{uy13}
Let $B^X: \Gamma (X,S(U,L)\otimes W)|_{\partial X} \longrightarrow
\Gamma (X,S(U,L)\otimes W)|_{\partial X}$
be the differential operator on $\partial X$ defined by
\begin{equation}\label{d2}
B^X=-\sum_{i=2}^{2l}c\left(\frac{\partial}{\partial
r}\right)c(e_i)\nabla_{e_i}^{S(U,L)\otimes W}.
\end{equation}
\end{defn}

By \cite{MR0397797},
$B^X$ is a formally self-adjoint first order elliptic differential
operator intrinsically defined on $\partial X$, which is the canonical boundary operator
associated to $D^X$ and preserves the natural
$\mathbb{Z}_2$-grading of $(S(U,L)\otimes W)|_{\partial X}$.

We now recall the Dirac type
operator \cite[Definition 1.1]{MR1870658} as well as
the boundary condition of Atiyah-Patodi-Singer type \cite{MR0397797}.

\begin{defn}\label{11y1}
By a Dirac type operator on $S(U,L)\otimes W$, we mean a
first order differential operator $D:\Gamma(X, S(U,L)\otimes W)\longrightarrow\Gamma(X, S(U,L)\otimes W)$ such that $D-D^X$ is an odd self-adjoint element of zeroth order, and that
its canonical boundary operator $B$ acting on
$\Gamma (X,S(U,L)\otimes W)|_{\partial X}$ is formally self-adjoint.
We will also call the restrictions $D_\pm$ of $D$ to $\Gamma (X,S_\pm(U,L)\otimes W)$ a
Dirac type operator.
\end{defn}

Let now $D$ be a $\mathbb{Z}/k$ Dirac type operator with its canonical boundary operator $B$. Obviously, $B$ preserves the
$\mathbb{Z}_2$-grading of $\Gamma (X,S(U,L)\otimes W)|_{\partial X}$,

Following \cite{MR0397797}, the boundary problem
\begin{equation}
\begin{split}
&(D_{+},P_{\geq 0,+}):\
\bigl\{s\,\big|\,s \in \Gamma(X,S_+(U,L)\otimes W), P_{\geq 0,+}(s\big|_{\partial X})=0 \bigr\}
\\
&\hspace{30pt}
\longrightarrow \Gamma(X,S_-(U,L)\otimes W),
\end{split}
\end{equation}
defines an elliptic boundary problem whose adjoint is $(D_-,P_{>0,-})$.
Moreover, it induces a Fredholm operator \cite{MR0397797}.
We will call the boundary problem $(D_{+},P_{\geq 0,+})$
the  Atiyah-Patodi-Singer boundary problem associated to $D_{+}$. Set
\begin{equation}\label{uy12}
\APS (D)=\dim \ker(D_{+},P_{\geq 0,+})-\dim \ker(D_{-},P_{>0,-}).
\end{equation}

\subsection{${\mathbb{Z}}/k$ circle actions and a variation formula}\label{sec2.2}
\begin{defn}
We will call a circle action on $X$ a $\mathbb{Z}/k$ circle action if
it preserves $\partial X$ and there exists a corresponding circle action on $Y$ such that
these two actions are compatible with $\pi$. The circle action is said to be nontrivial if it is not equal to identity.
\end{defn}

In what follows we assume that $X$ admits a nontrivial
$\mathbb{Z}/k$ circle action preserving the orientation and that the
$\mathbb{Z}/k$ circle action on $X$ lifts to $\mathbb{Z}/k$ circle
actions on $V$, $L$ and $W$, respectively. Without loss of generality, we may and
we will assume that these $\mathbb{Z}/k$ circle actions preserve
$g^{TX}$, $g^V$, $g^L$, $g^W$, $\nabla^V$, $\nabla^L$
and $\nabla^W$, respectively. We also assume that the $\mathbb{Z}/k$ circle actions on $TX$, $V$ and $L$ lift
to a $\mathbb{Z}/k$ circle action on $S(U,L)$ and preserves its $\mathbb{Z}_2$-grading.

Let $\mathscr{E}$ be a $\mathbb{Z}/k$ $S^1$-equivariant vector bundle over $X$.
Let $\mathscr{E}_Y$ be the $S^1$-equivariant vector bundle
over $Y$ induced from $\mathscr{E}$ through the map $\pi: \partial X \rightarrow Y$.
Recall
that the circle action on $\Gamma(X,\mathscr{E})$ is defined by $(g\cdot s)(x)=g(s(g^{-1}x))$ for
$g\in S^1$, $s\in \Gamma(X,\mathscr{E})$, $x\in X$.
Similarly, the group $S^1$ acts on $\Gamma(X,\mathscr{E})|_{\partial X}$ and $\Gamma(Y,\mathscr{E}_Y)$. For $\xi \in \mathbb{Z}$, by the weight-$\xi$ subspace
of $\Gamma(X,\mathscr{E})$ (resp. $\Gamma(X,\mathscr{E})|_{\partial X}$, $\Gamma(Y,\mathscr{E}_Y)$), we mean the subspace
of $\Gamma(X,\mathscr{E})$ (resp. $\Gamma(X,\mathscr{E})|_{\partial X}$, $\Gamma(Y,\mathscr{E}_Y)$) on which
$S^1$ acts as multiplication by $g^\xi$ for $g\in S^1$.

For any $\xi\in \mathbb{Z}$, let ${\rm E}^{\pm}_{\xi}$ (resp.
${\rm E}^{\pm}_{\xi,\partial}$, ${\rm E}^{\pm}_{Y,\xi}$) be the weight-$\xi$ subspaces of $\Gamma(X,S_{\pm}(U,L)\otimes W)$ (resp. $\Gamma(X,S_{\pm}(U,L)\otimes W)|_{\partial X}$, $\Gamma(Y, (S(U,L)\otimes W)_Y)$).

Let $D$ be a $\mathbb{Z}/k$ $S^1$-equivariant Dirac type operator on $\Gamma(S(U,L)\otimes W)$ with
canonical boundary operator $B$ acting on
$\Gamma(X,S(U,L)\otimes W)|_{\partial X}$. Let $P_{\geq0,+}$ be
the orthogonal projection associated to $B_+$.
For $\xi\in \mathbb{Z}$, let $D_{\pm,\xi}$ and $P_{\geq0,+,\xi}$ (resp. $P_{>0,-,\xi}$) be the restrictions of
$D_{\pm}$ and $P_{\geq0,+}$ (resp. $P_{>0,-}$) on the corresponding weight-$\xi$ subspaces ${\rm E}^{\pm}_{\xi}$
and ${\rm E}^{+}_{\xi,\partial}$ (resp. ${\rm E}^{-}_{\xi,\partial}$)
respectively. Then $(D_{+,\xi},P_{\geq 0,+,\xi})$
forms an elliptic boundary problem. Set
\begin{equation}\label{ly5}
\APS (D,\xi)=\dim \ker(D_{+,\xi},P_{\geq 0,+,\xi})-\dim \ker(D_{-,\xi},P_{>0,-,\xi})\ .
\end{equation}

Let $\bigl\{D_t: \Gamma(X,S(U,L)\otimes W)\longrightarrow \Gamma(X,S(U,L)\otimes W)\,\big|\,0 \leq t \leq1\bigr\}$
be a one parameter family of $\mathbb{Z}/k$ $S^1$-equivariant Dirac type operators with the canonical boundary operators $\bigl\{B_{t}\,\big|\, 0 \leq t \leq1 \bigr\}$. For any $t\in [0,1]$,
let $D^Y_{t,+}$ be the induced operator
from $B_{t,+}$ through the map $\pi: \partial X \rightarrow Y$, and let $B_{t,+,\xi}$ (resp. $D^Y_{t,+,\xi}$) be the restriction
of $B_{t,+}$ (resp. $D^Y_{t,+}$) on the weight-$\xi$
subspace ${\rm E}^{+}_{\xi,\partial}$ (resp. ${\rm E}^{+}_{Y,\xi}$).
We have the following variation formula.

\begin{thm}\label{thm1}{\rm (Compare with \cite[Theorem 1.2]{MR1870658})}
The following identity holds,
\begin{equation}
\begin{split}
\APS (D_1,\xi)-\APS (D_0,\xi)
&=-\SF\,\bigl\{B_{t,+,\xi}\,\big|\,0\leq t \leq1\bigr\}
\\
&=-k\SF\,\bigl\{D^Y_{t,+,\xi}\,\big|\,0\leq t \leq1\bigr\},
\end{split}
\end{equation}
where $\SF$ is the notation for the spectral flow of {\rm \cite{MR0397799}}.
In particular,
\begin{equation*}
\APS (D_1,\xi)\equiv\APS (D_0,\xi)\mod k\mathbb{Z}\,.
\end{equation*}
\end{thm}
\noindent{\bf Proof}\hspace{8pt} The proof is the same as that of \cite[Theorem 1.2]{MR1870658}.

\subsection{A mod $k$ localization formula for $\mathbb{Z}/k$ circle actions}\label{sec2.3}
Let $\mathcal{H}$ be the canonical basis of ${\rm Lie}(S^1)=\mathbb{R}$, i.e., for $t\in \mathbb{R}$,
$\exp(t\mathcal{H})=e^{2\pi \sqrt{-1}t}\in S^1$. Let $H$ be the Killing vector field on $X$
corresponding to $\mathcal{H}$. Since the circle action on $X$ is of $\mathbb{Z}/k$,
$H|_{\partial X}\subset T\partial X$ induces a Killing vector field
$H_Y$ on $Y$.
Let $X_H$ (resp. $Y_H$) be the zero set of $H$ (resp. $H_Y$) on $X$
(resp. $Y$). Then $X_H$ is a $\mathbb{Z}/k$ manifold and there is a
canonical map $\pi_{X_H}:\partial X_H\rightarrow Y_H$ induced by
$\pi$. In general, $X_H$ is not connected. We fix a connected
component $X_{H,\alpha}$ of $X_H$, and we omit the subscript
$\alpha$ if there is no confusion.

Clearly, $X_H$ intersects with $\partial X$ transversally. Let $g^{TX_H}$ be the metric on $X_H$ induced by $g^{TX}$.
Then $g^{TX_H}$ is naturally of product structure near $\partial X_{H}$. In fact, by choosing $\epsilon>0$ small enough, we know
$\mathscr{U}'_{\epsilon}=\mathscr{U}_{\epsilon}\cap X_H$ carries the metric naturally induced from
$g^{TX}|_{\mathscr{U}_{\epsilon}}$.

Let $\widetilde{\pi}:N\rightarrow X_H$ be the normal bundle to $X_H$ in $X$, which is identified
to be the orthogonal complement of $TX_{H}$ in $TX|_{X_H}$.
Then $TX|_{X_H}$ admits a $\mathbb{Z}/k$ $S^1$-equivariant decomposition
(cf. \cite[(1.8)]{MR2016198})
\begin{equation}\label{d3}
TX|_{X_H}=\bigoplus_{v\neq 0}N_v\oplus TX_{H}\,,
\end{equation}
where $N_v$ is a  $\mathbb{Z}/k$ complex vector bundle such that $g\in S^1$ acts on it by $g^v$ with $v\in \mathbb{Z}\backslash \{0\}$.
We will regard $N$ as a $\mathbb{Z}/k$ complex vector bundle and write
$N_\mathbb{R}$ for the underlying real vector bundle of $N$. Clearly, $N=\oplus_{v\neq 0}N_v$.
For $v\neq 0$, let $N_{v,\mathbb{R}}$ denote
the underlying real vector bundle of $N_v$.

Similarly, let
\begin{equation}\label{liu14}
W|_{X_H}=\bigoplus_{v}W_v, \quad V|_{X_H}=\bigoplus_{v\neq 0}V_v\oplus V_0^{\mathbb{R}}
\end{equation}
be the $\mathbb{Z}/k$ $S^1$-equivariant decompositions of the restrictions of $W$ and $V$ over $X_H$ respectively, where
$W_v$ and $V_v$ $(v\in \mathbb{Z})$ are $\mathbb{Z}/k$ complex vector bundles over $X_H$
on which $g\in S^1$ acts by $g^v$, and $V_0^{\mathbb{R}}$ is the real subbundle of $V$ such that
$S^1$ acts as identity. For $v\neq 0$, let $V_{v,\mathbb{R}}$ denote
the underlying real vector bundle of $V_v$. Denote by
$2p'=\dim V_0^{\mathbb{R}}$ and $2l'=\dim X_H$.

Let us write
\begin{align}\label{liu15}
L_{F}=L\otimes\left(\bigotimes_{v\neq0}\det
N_v\otimes\bigotimes_{v\neq0}\det V_v\right)^{-1}.
\end{align}
Then $TX_H\oplus V_0^{\R}$ has a $\mathbb{Z}/k$ Spin$^c$-structure since
$\omega_2(TX_H\oplus V_0^{\R})=c_1(L_F) \mod (2)$. Let $S(TX_H\oplus
V_0^{\R},L_F)$ be the fundamental spinor bundle for $(TX_H\oplus
V_0^{\R},L_F)$ as in Section \ref{sec2.1}.

Recall that $N_{v,\R}$ and $V_{v,\R}$ ($v\neq 0$) are canonically oriented by
their complex structures. The decompositions (\ref{d3}),
(\ref{liu14}) induce the orientations of $TX_H$ and $V_{0}^\mathbb{R}$ respectively. Let
$\{e_i\}_{i=1}^{2l'}$, $\{f_j\}_{j=1}^{2p'}$ be the corresponding
oriented orthonormal basis of $(TX_H,g^{TX_H})$ and
$(V_0^{\R},g^{V_0^{\R}})$. There are two canonical ways to consider
$S(TX_H\oplus V_0^{\R},L_F)$ as a $\mathbb{Z}_2$-graded vector
bundle . Let
\begin{align}\label{liu16}
\begin{split}
\tau_s&=(\sqrt{-1})^{l'}c(e_1)\cdots c(e_{2l'}),
\\
\tau_e&=(\sqrt{-1})^{l'+p'}c(e_1)\cdots c(e_{2l'})c(f_1)\cdots
c(f_{2p'})
\end{split}
\end{align}
be two involutions of $S(TX_H\oplus V_0^{\R},L_F)$. Then
$\tau_s^2=\tau_e^2=1$. We decompose $S(TX_H\oplus
V_0^{\R},L_F)=S_+(TX_H\oplus V_0^{\R},L_F)\oplus S_-(TX_H\oplus
V_0^{\R},L_F)$ corresponding to $\tau_s$ (resp. $\tau_e$) such that
$\tau_s|_{S_\pm(TX_H\oplus V_0^{\R},L_F)}=\pm 1$ (resp. $\tau_e|_{S_\pm(TX_H\oplus V_0^{\R},L_F)}=\pm 1$).

Let $C(N_\mathbb{R})$ be the Clifford algebra bundle of $(N_\mathbb{R},g^N)$.
Then $\Lambda(\overline{N}^*)$ is a $C(N_\mathbb{R})$-Clifford module.
Namely, for $e\in N$, let $e'\in \overline{N}^*$ correspond to $e$ by
the metric $g^N$, and let
\begin{equation}\label{liu19}
c(e)=\sqrt{2}\,e'\wedge\ ,\quad c(\overline{e})=-\sqrt{2}\,i_{\overline{e}}\ ,
\end{equation}
where $\wedge$ and $i$ denote the exterior and interior multiplications, respectively.
Let $\tau^N$ be the involution on $\Lambda(\overline{N}^*)$ given by
$\tau^N|_{\Lambda^{{\rm even}/{\rm odd}}(\overline{N}^*)}=\pm 1$.

Similarly, we can define the Clifford action of $C(V_{v,\mathbb{R}})$ on the $C(V_{v,\mathbb{R}})$-Clifford module
$\Lambda(\overline{V}_v^*)$ with the involution
$\tau^V_v|_{\Lambda^{{\rm even}/{\rm odd}}(\overline{V}_v^*)}=\pm 1$.

Upon restriction to $X_H$, one has the following $\mathbb{Z}/k$ isomorphisms of
$\mathbb{Z}_2$-graded Clifford modules over $X_H$ (compare with \cite[(1.49)]{MR2016198}),
\begin{equation}\label{liu17}
\bigl(S(U,L),\tau_s\bigr)|_{X_H}\simeq \bigl(S(TX_H\oplus
V_0^{\R},L_F),\tau_s\bigr)\,\widehat{\otimes}\,\bigl(\Lambda\overline{N}^*,\tau^N\bigr)
\widehat{\otimes}\,\widehat{\bigotimes_{v\neq0}}\,\bigl(\Lambda \overline{V}^*_v,{\rm id}\bigr),
\end{equation}
where ${\rm id}$ denotes the trivial involution, and
\begin{equation}\label{liu18}
\bigl(S(U,L),\tau_e\bigr)|_{X_H}\simeq \bigl(S(TX_H\oplus
V_0^{\R},L_F),\tau_e\bigr)\,\widehat{\otimes}\,\bigl(\Lambda\overline{N}^*,\tau^N\bigr)
\widehat{\otimes}\,\widehat{\bigotimes_{v\neq0}}\,\bigl(\Lambda \overline{V}^*_v,\tau_v^V\bigr).
\end{equation}
Here we denote by
$\widehat{\otimes}$ the $\mathbb{Z}_2$-graded tensor product  (cf. \cite[pp. 11]{MR1031992}). Furthermore,
isomorphisms (\ref{liu17}), (\ref{liu18}) give the identifications of the canonical connections on
the bundles (compare with \cite[(1.13)]{MR2016198}). We still denote the involution on $\bigl(S(TX_H\oplus
V_0^{\R},L_F)$ by $\tau$.

Let $R$ be a $\mathbb{Z}/k$ Hermitian vector bundle over $X_H$ endowed with a $\mathbb{Z}/k$
Hermitian connection. We make the assumption that the Hermitian metric and the Hermitian connection
are both of product structure near $\partial X_H$. We will denote by $D^{X_H}\otimes R$ the twisted Spin$^c$
Dirac operator on $S(TX_H\oplus V_0^{\R},L_F)\otimes R$ and by $D^{X_{H,\alpha}}\otimes R$ its restriction
to $X_{H,\alpha}$ (cf. Definition \ref{diracdef}).

We denote by $K(X_{H})$ the $K$-group of $\mathbb{Z}/k$
complex vector bundles over $X_{H}$ (cf. \cite[pp. 285]{MR1144425}).
We use the same notations as in \cite[pp. 128]{MR2016198},
\begin{equation}\label{liu27}
\begin{split}
{\rm Sym}_q(R)&=\sum_{n=0}^{+\infty}q^n{\rm Sym}^n(R)\in K(X_H)[[q]],
\\
\Lambda_q(R)&=\sum_{n=0}^{+\infty}q^n\Lambda^n(R)\in K(X_H)[[q]],
\end{split}
\end{equation}
for the symmetric and exterior power
operations in $K(X_H)[[q]]$, respectively.

Let $S^1$ act on $L|_{X_H}$ by sending $g\in S^1$ to $g^{l_c}$ ($l_c\in
\mathbb{Z}$) on $X_H$. Then $l_c$ is locally constant on $X_H$.
Following \cite[(1.50)]{MR2016198}, we define the following elements in $K(X_H)[[q^{\frac{1}{2}}]]$,
\begin{equation}\label{liu25}
\begin{split}
R_{\pm}(q)&=q^{\tfrac{1}{2}\sum_v|v|\dim N_v-\tfrac{1}{2}\sum_v v\dim V_v+\tfrac{1}{2}l_c}
\bigotimes_{v>0}\left({\rm
Sym}_{q^v}(N_v)\otimes\det N_v\right)
\\
&\quad \otimes\bigotimes_{v<0}{\rm
Sym}_{q^{-v}}(\overline{N}_v)\otimes\bigotimes_{v\neq 0}\Lambda_{\pm q^v}(V_v)\otimes\Bigl(\,\sum_vq^vW_v\Bigr)
\\
&=\sum_{n} R_{\pm,n}q^n\ ,
\end{split}
\end{equation}
\begin{equation}\label{liu26}
\begin{split}
R'_{\pm}(q)&=q^{-\tfrac{1}{2}\sum_v|v|\dim N_v-\tfrac{1}{2}\sum_v v\dim V_v+\tfrac{1}{2}l_c}
\bigotimes_{v>0}{\rm
Sym}_{q^{-v}}(\overline{N}_v)
\\
&\quad \otimes\bigotimes_{v<0}\left({\rm
Sym}_{q^v}(N_v)\otimes\det N_v\right)\otimes\bigotimes_{v\neq 0}\Lambda_{\pm q^v}(V_v)\otimes\Bigl(\,\sum_vq^vW_v\Bigr)
\\
&=\sum_{n} R'_{\pm,n}q^n.
\end{split}
\end{equation}
As explained in \cite[pp. 139]{MR2016198}, since $TX\oplus V\oplus L$ is spin,
one gets
\begin{equation}
\sum_vv\dim N_v+\sum_v v\dim V_v+l_c \equiv 0 \mod (2).
\end{equation}
Therefore, $R_{\pm,\xi}(q)$, $R'_{\pm,\xi}(q)\in K(X_H)[[q]]$.

Clearly each $R_{\pm,\xi}$, $R'_{\pm,\xi}$ $(\xi\in \mathbb{Z})$ is a $\mathbb{Z}/k$ vector
bundle over $X_H$ carrying a canonically induced $\mathbb{Z}/k$
Hermitian metric and a canonically induced $\mathbb{Z}/k$ Hermitian connection,
which are both of product structure near $\partial X_H$.

We now state a mod $k$ localization formula which generalizes \cite[Theorem 1.2]{MR2016198}
to the case of $\mathbb{Z}/k$ manifolds. It also generalizes the $\mathbb{Z}/k$ equivariant index theorem in
\cite[Theorem 2.1]{MR1993996} to the case of Spin$^c$-manifolds.

\begin{thm}\label{liu28}
For any $\xi \in\mathbb{Z}$, the following identities hold,
\begin{align}\label{liu29}
&\APS_{\tau_s}\bigl(D^X,\xi\bigr)
\equiv\sum_{\alpha}(-1)^{\sum_{0<v}\dim N_v}\APS_{\tau_s}\bigl(D^{X_{H,\alpha}}\otimes R_{+,\xi}\bigr)
\notag\\
&\hspace{30pt}
\equiv\sum_{\alpha}(-1)^{\sum_{v<0}{\rm
dim}\,N_v}\APS_{\tau_s}\bigl(D^{X_{H,\alpha}}\otimes R'_{+,\xi}\bigr)
\mod k\mathbb{Z}\ ,
\\
\label{liu30}
&\APS_{\tau_e}\bigl(D^X,\xi\bigr)
\equiv\sum_{\alpha}(-1)^{\sum_{0<v}\dim N_v}\APS_{\tau_e}\bigl(D^{X_{H,\alpha}}\otimes R_{-,\xi}\bigr)
\notag\\
&\hspace{30pt}
\equiv\sum_{\alpha}(-1)^{\sum_{v<0}{\rm
dim}\,N_v}\APS_{\tau_e}\bigl(D^{X_{H,\alpha}}\otimes R'_{-,\xi}\bigr)
\mod k\mathbb{Z}\ .
\end{align}
\end{thm}
\noindent{\bf Proof}\hspace{8pt} The proof will be given in Section 3.

\subsection{A $\mathbb{Z}/k$ extension of Hattori's vanishing theorem}\label{sec2.4}

In this subsection, we assume that $TX$ has a $\mathbb{Z}/k$ $S^1$-equivariant  almost complex structure $J$.
Then one has the canonical splitting
\begin{equation}
TX\otimes_{\mathbb{R}}\mathbb{C}=T^{(1,0)}X\oplus T^{(0,1)}X,
\end{equation}
where $T^{(1,0)}X$ and $T^{(0,1)}X$ are the eigenbundles of $J$ corresponding to the eigenvalues $\sqrt{-1}$
and $-\sqrt{-1}$, respectively.

Let $K_{X}=\det(T^{(1,0)}X)$ be the determinant line bundle of $T^{(1,0)}X$ over $X$. Then the
complex spinor bundle $S(TX,K_{X})$ for $(TX,K_{X})$ is $\Lambda(T^{*(0,1)}X)$ (cf. \cite[Appendix D]{MR1031992}).

We suppose that $c_1(T^{(1,0)}X)=0 \mod (N)$ $(N\in\mathbb{Z},N\geq 2)$.
As explained in Section \ref{sec2.1}, the complex line bundle $K_{X}^{1/N}$ is
well defined over $X$. After replacing the $S^1$ action by its $N$-fold action, we can always assume that $S^1$ acts on
$K_{X}^{{1/ N}}$. For $s\in \mathbb{Z}$, let $D^X\otimes K_{X}^{s/N}$ be the twisted Spin$^c$ Dirac operator
on $\Lambda(T^{*(0,1)}X)\otimes K_{X}^{{s/N}}$ defined as in (\ref{d1}).

Using Theorem \ref{liu28}, we can generalize the main result of Hattori \cite{MR508087} to the case of $\mathbb{Z}/k$ almost complex manifolds.

\begin{thm}\label{main}
Assume that $X$ is a connected  $\mathbb{Z}/k$ almost complex manifold
with a nontrivial $\mathbb{Z}/k$ circle action.
If $c_1(T^{(1,0)}X)=0 \mod (N)$ $(N\in\mathbb{Z},N\geq 2)$, then for $s\in \mathbb{Z}$, $-N<s<0$,
$D^X\otimes K_{X}^{s/N}$ has vanishing property in $\mathbb{Z}/k$ category. In particular, the following identity holds,
\begin{equation}\label{uy1}
\topo\bigl(D^X\otimes K_{X}^{s/N}\bigr)= 0 \quad \text{in}\quad \mathbb{Z}/k\mathbb{Z}\ .
\end{equation}
\end{thm}

\noindent{\bf Proof}\hspace{8pt}
Using the almost complex structure on $TX_{H}$ induced by the
almost complex structure $J$ on $TX$ and by (\ref{d3}), we know
\begin{equation}
T^{(1,0)}X\big|_{X_H}=\bigoplus_{v\neq 0}N_v\oplus T^{(1,0)}X_{H}\,,
\end{equation}
where $N_v$ are complex subbundles of $T^{(1,0)}X\big|_{X_H}$ on which $g\in S^1$ acts by multiplication
by $g^v$.

We claim that for each $\xi\in \mathbb{Z}$, the following identity holds,
\begin{equation}\label{ly2}
\APS \bigl(D^X\otimes K_{X}^{s/N},\xi\bigr)\equiv0\ \mod k\mathbb{Z}\ .
\end{equation}

In fact, if $X_H=\emptyset$, the empty set, by Theorem \ref{liu28}, (\ref{ly2}) is
obvious.

When $X_H\neq\emptyset$, we see that $\sum_v|v|\dim N_v>0$ (i.e., at least one of the $N_v$'s is nonzero) on
each connected component of $X_H$. Consider $R_{+}(q)$, $R'_{+}(q)$ of (\ref{liu25}) and (\ref{liu26}) for the case
that $V=0$, $W=K_{X}^{s/N}$. We deduce that
\begin{equation*}
\begin{split}
R_{+,\xi}=0 &\,\quad\text{if}\quad\xi<a_1=\inf_\alpha\Bigl(\,\frac{1}{2}\sum_v|v|\dim N_v
+\Bigl(\frac{1}{2}+\frac{s}{N}\Bigr)\sum_vv\dim N_v\Bigr),
\\
R'_{+,\xi}=0 &\,\quad\text{if}\quad\xi>a_2=\sup_\alpha\Bigl(-\frac{1}{2}\sum_v|v|\dim N_v
+\Bigl(\frac{1}{2}+\frac{s}{N}\Bigr)\sum_vv\dim N_v\Bigr).
\end{split}
\end{equation*}
Since $-N<s<0$, we know $a_1>0$ and $a_2<0$.
By using Theorem \ref{liu28}, we see that (\ref{ly2}) holds for any $\xi\in\mathbb{Z}$.

Now Theorem \ref{main} follows easily from \eqref{fm1}, \eqref{fm2} and \eqref{ly2}.
The proof of Theorem \ref{main} is completed.

\begin{rem}
From the proof of Theorem \ref{main}, one also deduces that if $X$ is a connected
$\mathbb{Z}/k$ almost complex manifold with a nontrivial $\mathbb{Z}/k$ circle action, then
$D^X$, $D^X\otimes K^{-1}_{X}$ are rigid in $\mathbb{Z}/k$ category.
\end{rem}

\section{A proof of Theorem \ref{liu28}}\label{sec3}
In this section, following Zhang \cite{MR1993996} and by making use of the analysis of Wu-Zhang \cite{MR1601858}
and Dai-Zhang \cite{MR1870658} as well as Liu-Ma-Zhang \cite{MR2016198}, which
in turn depend on the analytic localization techniques of Bismut-Lebeau \cite{MR1188532}, we present a proof of Theorem \ref{liu28}.

This section is organized as follows.
In Section \ref{sec3.1}, we recall a result from \cite{MR1601858} concerning the Witten deformation on flat spaces.
In Section \ref{sec3.2}, we establish the Taylor expansions of $D^X$ and $c(H)$
(resp. $B^X$) near the fixed point set $X_H$ (resp. $\partial X_H$).
In Section \ref{sec3.3}, following \cite[Section 3(b)]{MR1870658},
we decompose the Dirac type operators under consideration to a sum of four operators and introduce a deformation of the Dirac
type operators as well as their associated boundary operators.
In Section \ref{sec3.4}, by using the techniques of \cite[Section 3(c)]{MR1870658},
\cite[Section 1.2]{MR2016198} and \cite[Section 9]{MR1188532},
we carry out various estimates for certain operators and prove the Fredholm property of the Atiyah-Patodi-Singer
type boundary problem for the deformed operators introduced in Section \ref{sec3.3}.
In Section \ref{sec3.5}, we complete the proof of Theorem \ref{liu28}.

\subsection{Witten deformation on flat spaces}\label{sec3.1}

Recall that $\mathcal{H}$ is the canonical basis of ${\rm Lie}(S^1)=\mathbb{R}$.
In this subsection, let $W$ be a complex vector space of dimension $n$ with an Hermitian form. Let
$\rho$ be a unitary representation of the circle group $S^1$ on $W$
such that all the weights are nonzero. Suppose $W^{\pm}$ are the
subspaces of $W$ corresponding to the positive and negative weights
respectively, with $\dim_{\mathbb{C}}W^-=\nu$, $\dim_{\mathbb{C}}W^+=n-\nu$. Let $z=\{z^i\}$ be the complex linear
coordinates on $W$ such that the Hermitian structure on $W$ takes
the standard form and $\rho$ is diagonal with weights
$\lambda_i\in\mathbb{Z}\backslash\{0\}\,(1\leq i\leq n)$, and
$\lambda_i<0$ for $i\leq \nu$. The Lie algebra action on $W$ is
given by the vector field
\begin{equation}\label{liu31}
H=2\pi\sqrt{-1}\sum_{i=1}^n
\lambda_i\left(z^i\frac{\partial}{\partial
z^i}-\bar{z}^i\frac{\partial}{\partial \bar{z}^i}\right).
\end{equation}
Set
\begin{equation}\label{liu32}
K^{\pm}(W)={\rm Sym}((W^{\pm})^*)\otimes{\rm
Sym}(W^{\mp})\otimes\det(W^{\mp}).
\end{equation}
Let $E$ be a finite dimensional complex vector space with an
Hermitian form and suppose $E$ carries a unitary representation of
$S^1$.

Let $\overline{\partial}$ be the twisted Dolbeault operator acting on
$\Omega^{0,*}(W,E)$, the set of smooth sections of
$\Lambda(\overline{W}^*)\otimes E$ on $W$. Let $\overline{\partial}^*$ be
the formal adjoint of $\overline{\partial}$. Let $D=\sqrt{2}(\overline{\partial}+\overline{\partial}^*)$.
Let $c(H)$ be the Clifford action of $H$ on
$\Lambda(\overline{W}^*)$ defined as in  (\ref{liu19}). Let $\mathscr{L}_H$ be the Lie
derivative along $H$ acting on $\Omega^{0,*}(W,E)$.

The following result was proved in \cite[Proposition 3.2]{MR1601858}.
\begin{prop}\label{liu33}
{\rm 1}. A basis of the space of $L^2$-solutions of $D+\sqrt{-1}c(H)$
(resp. $D-\sqrt{-1}c(H)$) on the space of $C^{\infty}$ sections of
$\Lambda(\overline{W}^*)$ is given by
\begin{equation}\label{liu34}
\big(\prod_{i=1}^{\nu}z_i^{k_i}\big)\big(\prod_{i=\nu+1}^{n}
\bar{z}_i^{k_i}\big)e^{-\sum_{i=1}^n\pi|\lambda_i||z_i|^2}
d\bar{z}_{\nu+1}\cdots d\bar{z}_n\quad (k_i\in\mathbb{N})
\end{equation}
with weight $\sum_{i=1}^{\nu}k_i|\lambda_i|+\sum_{i=\nu+1}^{n}(k_i+1)|\lambda_i|$
(resp.
\begin{equation}\label{liu35}
\big(\prod_{i=1}^{\nu}\bar{z}_i^{k_i}\big)\big(\prod_{i=\nu+1}^{n}
z_i^{k_i}\big)e^{-\sum_{i=1}^n\pi|\lambda_i||z_i|^2}
d\bar{z}_{1}\cdots d\bar{z}_{\nu}\quad (k_i\in\mathbb{N})
\end{equation}
with weight
$-\sum_{i=\nu+1}^{n}k_i|\lambda_i|-\sum_{i=1}^{\nu}(k_i+1)|\lambda_i|$).

So the space of $L^2$-solution of a given weight of
$D+\sqrt{-1}c(H)$ (resp. $D-\sqrt{-1}c(H)$) on the space of
$C^{\infty}$ sections of $\Lambda(\overline{W}^*)\otimes E$ is finite
dimensional. The direct sum of these weight spaces is isomorphic to
$K^-(W)\otimes E$ (resp. $K^+(W)\otimes E$) as representations of
$S^1$.

{\rm 2}. When restricted to an eigenspace of $\mathscr{L}_H$, the operator
$D+\sqrt{-1}c(H)$ (resp. $D-\sqrt{-1}c(H)$) has discrete
eigenvalues.

\end{prop}

\subsection{A Taylor expansion of certain operators near the fixed-point set}\label{sec3.2}

Following \cite[Section 8(e)]{MR1188532}, we now describe a coordinate system on $X$ near $X_H$.
For $\varepsilon>0$, set ${\mathscr B}_{\varepsilon}=\bigl\{Z\in N\big|\ |Z|<\varepsilon\bigr\}$.
Since $X$ and $X_H$ are compact, there exists $\varepsilon_0>0$ such
that for $0<\varepsilon\leq\varepsilon_0$, the exponential map
\begin{equation*}
(y,Z)\in N\longmapsto \exp_y^X(Z)\in X
\end{equation*}
is a diffeomorphism from
${\mathscr B}_{\varepsilon}$ onto a tubular neighborhood ${\mathscr V}_{\varepsilon}$ of
$X_H$ in $X$. From now on, we identify ${\mathscr B}_{\varepsilon}$ with
${\mathscr V}_{\varepsilon}$ and use the notation $x=(y,Z)$ instead of
$x=\exp_y^X(Z)$. Finally, we identify $y\in X_H$ with $(y,0)\in N$.

Let $\widetilde{\pi}^*((S(U,L)\otimes W)|_{X_H})$ be the vector bundle on
$N$ obtained by pulling back $(S(U,L)\otimes W)|_{X_H}$ for
$\widetilde{\pi}:N\rightarrow X_H$.

Let $g^{TX_H}$, $g^{N}$ be the corresponding metrics on $TX_H$ and
$N$ induced by $g^{TX}$. Let ${\rm d}v_X$, ${\rm d}v_{X_H}$ and ${\rm d}v_N$ be the
corresponding volume elements on $(TX,g^{TX})$, $(TX_H,g^{TX_H})$
and $(N,g^{N})$. Let $k(y,Z)$ ($(y,Z)\in {\mathscr B}_{\varepsilon}$) be the
smooth positive function defined by
\begin{equation}\label{liu38}
{\rm d}v_X(y,Z)=k(y,Z){\rm d}v_{X_H}(y){\rm d}v_{N_y}(Z).
\end{equation}
Then $k(y)=1$ and $\frac{\partial k}{\partial Z}(y)=0$ for $y\in
X_H$. The latter follows from the well-known fact that $X_H$ is totally
geodesic in $X$.

For $x=(y,Z)\in {\mathscr V}_{\varepsilon_0}$, we will identify $S(U,L)_x$ with
$S(U,L)_y$ and $W_x$ with $W_y$ by the parallel transport with
respect to the $S^1$-invariant connections $\nabla^{S(U,L)}$ and
$\nabla^W$ respectively, along the geodesic $t\longmapsto (y,tZ)$.
The induced identification of $(S(U,L)\otimes W)_x$ with
$(S(U,L)\otimes W)_y$ preserves the metric and the
$\mathbb{Z}_2$-grading, and moreover, is $S^1$-equivariant.
Consequently, $D^X$ can be considered as an operator acting on the
sections of the bundle $\widetilde{\pi}^*((S(U,L)\otimes W)|_{X_H})$ over
${\mathscr B}_{\varepsilon_0}$ commuting with the circle action.

For $\varepsilon>0$, let $\textbf{E}({\varepsilon})$ (resp.
$\textbf{E}$) be the set of smooth sections of
$\widetilde{\pi}^*((S(U,L)\otimes W)|_{X_H})$ on ${\mathscr B}_{\varepsilon}$ (resp.
on the total space of $N$). If $f,g\in\textbf{E}$ have compact
supports, we will write
\begin{equation}\label{liu39}
\langle f,g\rangle
=\left(\frac{1}{2\pi}\right)^{\dim X}\int_{X_H}\left(\int_N\langle
f,g\rangle(y,Z){\rm d}v_{N_y}(Z)\right){\rm d}v_{X_H}(y).
\end{equation}
Then $k^{1/2}D^Xk^{-1/2}$ is a (formally) self-adjoint operator on
$\textbf{E}$.

The connection $\nabla^N$ on $N$ induces a splitting $TN=N\oplus
T^HN$, where $T^HN$ is the horizontal part of $TN$ with respect to
$\nabla^N$. Moreover, since $X_H$ is totally geodesic, this
splitting, when restricted to $X_H$, is preserved by the connection
$\nabla^{TX}$ on $TX|_{X_H}$. Let $\widetilde{\nabla}$ be the connection
on $(S(U,L)\otimes W)|_{X_H}$ induced by the restriction of
$\nabla^{S(U,L)\otimes W}$ to $X_H$. We denote by
$\widetilde\pi^*\widetilde{\nabla}$ the pulling back  of the connection $\widetilde{\nabla}$ on
$(S(U,L)\otimes W)|_{X_H}$ to the bundle
$\widetilde{\pi}^*((S(U,L)\otimes W)|_{X_H})$.

We choose a local orthonormal basis of $TX$ such that
$e_1,\cdots,e_{2l'}$ form a basis of $TX_H$, and
$e_{2l'+1},\cdots,e_{2l}$, that of $N_{\mathbb{R}}$. Denote the horizontal lift
of $e_{i}\,(1\leq i\leq 2l')$ to $T^HN$ by $e^H_i$. We define
\begin{equation}\label{liu40}
D^H=\sum_{i=1}^{2l'}c(e_i)(\widetilde\pi^*\widetilde{\nabla})_{e^H_i}, \quad
D^N=\sum_{i=2l'+1}^{2l}c(e_i)(\widetilde\pi^*\widetilde{\nabla})_{e_i}.
\end{equation}
Clearly, $D^N$ acts along the fibers of $N$. Let $\overline{\partial}^N$
be the $\overline{\partial}$-operator along the fibers of $N$, and let
$\overline{\partial}^{N*}$ be its formal adjoint with respect to
(\ref{liu39}). It is easy to see that
$D^N=\sqrt{2}\bigl(\overline{\partial}^N+\overline{\partial}^{N*}\bigr)$. Both $D^N$ and
$D^H$ are formally self-adjoint with respect to (\ref{liu39}).

For $T>0$, we define a scaling
$f\in\textbf{E}({\varepsilon_0})\rightarrow
S_Tf\in\textbf{E}({\varepsilon_0\sqrt{T}})$ by
\begin{equation}\label{liu41}
S_Tf(y,Z)=f\left(y,\frac{Z}{\sqrt{T}}\right),\quad (y,Z)\in
{\mathscr B}_{\varepsilon_0\sqrt{T}}\ \ .
\end{equation}
For a first order differential operator
\begin{equation}\label{liu42}
Q_T=\sum_{i=1}^{2l'}a^i_T(y,Z)(\widetilde\pi^*\widetilde{\nabla})_{e^H_i}+\sum_{i=2l'+1}^{2l}
b^i_T(y,Z)(\widetilde\pi^*\widetilde{\nabla})_{e_i}+c_T(y,Z)
\end{equation}
acting on $\textbf{E}({\varepsilon_0\sqrt{T}})$, where $a_T^i,
b_T^i$, and $c_T$ are endomorphisms of $\widetilde{\pi}^*((S(U,L)\otimes
W)|_{X_H})$ which depend smoothly on $(y,Z)$, we write
\begin{equation}\label{liu43}
Q_T=O\bigl(|Z|^2\partial^N+|Z|\partial^H+|Z|+|Z|^p\bigr),
\end{equation}
if there is a constant $C>0$, $p\in
\mathbb{N}$ such that for any $T\geq 1$, $(y,Z)\in
{\mathscr B}_{\varepsilon_0\sqrt{T}}$, we have
\begin{equation}
\begin{split}
|a_T^i(y,Z)|&\leq C|Z|\quad (1\leq i\leq 2l'),
\\
|b_T^i(y,Z)|&\leq C|Z|^2\quad (2l'+1\leq i\leq 2l),
\\
|c_T(y,Z)|&\leq C(|Z|+|Z|^p)\,.
\end{split}
\end{equation}

Let $\mathbf{E}_\partial$ be the set of smooth sections of
$\widetilde{\pi}^*((S(U,L)\otimes W)|_{X_H})$ over $N|_{\partial X_H}$.
On the boundary of $X_H$, we choose the local orthonormal basis as in Definition \ref{uy13}.
Similarly as in \eqref{d2}, we define
\begin{equation}
B^H=-\sum_{i=2}^{2l'}c\left(\frac{\partial}{\partial r}\right)c(e_i)(\widetilde\pi^*\widetilde{\nabla})_{e^H_i},
\quad
B^N=-c\left(\frac{\partial}{\partial r}\right)D^N|_{\partial X_H}
\end{equation}
on $\mathbf{E}_\partial$ (compare with \eqref{liu40}).

Let $J_H$ be the representation of Lie$(S^1)$ on $N$. Then
$Z\rightarrow J_HZ$ is a Killing vector field on $N$.
We have the following analogue of \cite[Theorem 8.18]{MR1188532}, \cite[Proposition 1.2]{MR2016198} and \cite[Proposition 3.3]{MR1601858}.
\begin{prop}\label{liu45}
As $T\rightarrow +\infty$,
\begin{equation*}
\begin{split}
S_Tk^{1/2}D^Xk^{-1/2}S_T^{-1}
&=\sqrt{T}D^N+D^H+\frac{1}{\sqrt{T}}
O(|Z|^2\partial^N+|Z|\partial^H+|Z|),
\\
S_Tk^{1/2}c(H)k^{-1/2}S_T^{-1}
&=\frac{1}{\sqrt{T}}c(J_HZ)+\frac{1}{\sqrt{T^3}}O(|Z|^3),
\\
S_Tk^{1/2}B^Xk^{-1/2}S_T^{-1}
&=\sqrt{T}B^N+B^H+\frac{1}{\sqrt{T}}
O(|Z|^2\partial^N+|Z|\partial^H+|Z|).
\end{split}
\end{equation*}
\end{prop}

\subsection{A decomposition of Dirac type operators under consideration and the associated deformation}\label{sec3.3}

For $p\geq 0$, let ${\rm E}^p$ (resp. ${\rm E}_\partial^p$, $\textbf{E}^p$, ${\rm F}^p$, ${\rm F}^p_{\partial}$) be the set
of sections of the bundles $S(U,L)\otimes W$ over $X$ (resp. $(S(U,L)\otimes W)|_{\partial X}$ over $\partial X$,
$\widetilde{\pi}^*((S(U,L)\otimes W)|_{X_H})$ over $N$, $S(TX_H\oplus V^{\mathbb{R}}_0,L_F)\otimes K^-(N)\otimes
\bigl(\widehat{\otimes}_{v\neq 0}\Lambda V_v\otimes W\bigr)|_{X_H}$ over $X_H$, $\bigl(S(TX_H\oplus V^{\mathbb{R}}_0,L_F)\otimes K^-(N)\otimes
\widehat{\otimes}_{v\neq 0}\Lambda V_v\otimes W\bigr)|_{\partial X_H}$ over $\partial X_H$) which lie in the $p$-th
Sobolev spaces. The group $S^1$ acts on all these spaces (cf. Section \ref{sec2.2}). For
any $\xi\in\mathbb{Z}$, let ${\rm E}^p_{\xi}$, ${\rm E}_{\xi,\partial}^p$ ${\bf E}^p_{\xi}$,
${\rm F}^p_{\xi}$ and ${\rm F}^p_{\xi,\partial}$ be the corresponding weight-$\xi$ subspaces, respectively.

Recall that the constant $\varepsilon_0>0$ is
defined in last subsection. We now take
$\varepsilon\in(0,\frac{\varepsilon_0}{2}]$. Let $\rho:\R\rightarrow[0,1]$ be a smooth
function such that
\begin{equation}\label{liu46}
\rho(a)=
\begin{cases}1 & \text{if $a\leq \frac{1}{2}$,}\\0 & \text{if $a\geq 1$.}
\end{cases}
\end{equation}
For $Z\in N$, set $\rho_\varepsilon(Z)=\rho(\frac{|Z|}{\varepsilon})$.

By Proposition \ref{liu33}, the solution space of the operator $D^N+\sqrt{-1}Tc(J_HZ)$
along the fiber $N_y$ ($y\in X_{H}$) is the $L^2$ completion of $K^-(N_y)\otimes(\widehat{\otimes}_{v\neq 0}\Lambda V_v\otimes W)_y$.
They form an infinite dimensional Hermitian complex vector bundle $K^-(N)\otimes(\widehat{\otimes}_{v\neq 0}\Lambda V_v\otimes W)|_{X_H}$ over $X_H$, with the Hermitian connection induced from those on $N$, $V|_{X_H}\rightarrow X_H$ and $W|_{X_H}\rightarrow X_H$. Let $\theta$ be the isomorphism from $L^2(X_H,K^-(N)\otimes(\widehat{\otimes}_{v\neq 0}\Lambda V_v\otimes W)|_{X_H})$ to $L^2(N,\widetilde\pi^*((\Lambda \overline{N}^*\otimes\widehat{\otimes}_{v\neq 0}\Lambda V_v\otimes W)|_{X_H}))$ given by Proposition \ref{liu33}.

Let $\alpha\in \Gamma\left(X_H,S(TX_H\oplus V^{\mathbb{R}}_0,L_F)\right)$, $\phi\in L^2(X_H,K^-(N)\otimes(\widehat{\otimes}_{v\neq 0}\Lambda V_v\otimes W)|_{X_H})$,
$\sigma=\alpha\otimes\phi$. We define a linear map
\begin{equation}\label{itxi}
I_{T,\xi}:{\rm F}^p_{\xi}\longrightarrow{\rm\bf E}^p_{\xi},
\qquad\sigma\longmapsto T^{\frac{\dim N_\R}{2}}\rho_{\varepsilon}(Z)\,\widetilde{\pi}^*\alpha\wedge
S_T^{-1}(\theta\phi).
\end{equation}
In general, there exist
$c(\varepsilon )>0$ and $C>0$ such that $c(\varepsilon )<\|I_{T,\xi}\|<C$.

Let the image of $I_{T,\xi}$ from  ${\rm F}_{\xi}^p$ be
${\rm \bf E}_{T,\xi}^p=I_{T,\xi}{\rm F}_{\xi}^p\subseteq {\rm\bf E}^p_\xi$. Denote the orthogonal complement of
${\bf E}_{T,\xi}^0$ in ${\bf E}_\xi^0$ by
${\bf E}_{T,\xi}^{0,\bot}$, and let
${\bf E}_{T,\xi}^{p,\bot}=\textbf{E}_{\xi}^{p}\cap\textbf{E}_{T,\xi}^{0,\bot}$.
Let $p_{T,\xi}$ and $p_{T,\xi}^\bot$ be the orthogonal projections from
$\textbf{E}_\xi^0$ to $\textbf{E}_{T,\xi}^0$ and $\textbf{E}_{T,\xi}^{0,\bot}$ respectively.

We denote by $\bigl(\bigl(\,\widehat{\bigotimes}_{v\neq 0}\Lambda V_v\,\bigr)\otimes
\bigl(\,\bigoplus_v W_v\,\bigr)\bigr)_{\xi-\frac{1}{2}\sum_{v}|v|\dim N_{v}}$ the subbundle of
$\bigl(\,\widehat{\bigotimes}_{v\neq 0}\Lambda V_v\,\bigr)\otimes\bigl(\,\bigoplus_v W_v\,\bigr)$ whose weight equals to
$\xi-\tfrac{1}{2}\sum_{v}|v|\dim N_{v}$ with respect to the given circle action.
Let $q_{\xi}$ be the orthogonal bundle projection from the vector bundle
\begin{equation*}
\bigl(\,\widehat{\bigotimes}_{v\neq 0}\Lambda \overline{N}^*_v\,\bigr)\,\widehat\otimes\,
\bigl(\,\widehat{\bigotimes}_{v\neq 0}\Lambda V_v\,\bigr)\otimes \bigl(\,\bigoplus_v W_v\,\bigr)\longrightarrow X_H
\end{equation*}
to its subbundle
\begin{equation*}
\bigotimes_{v>0}\det N_{v}\,\widehat\otimes\,
\Bigl(\bigl(\,\widehat{\bigotimes}_{v\neq 0}\Lambda
V_v\,\bigr)\otimes \bigl(\,\bigoplus_v
W_v\,\bigr)\Bigr)_{\xi-\frac{1}{2}\sum_{v}|v|\dim
N_{v}}\longrightarrow X_H\, .
\end{equation*}

We now proceed to deduce a formula which computes $p_{T,\xi}s$ for $s\in{\rm E}^0_{\xi}$ explicitly under a
local unitary trivialization of $N$.

For $y_0\in X_{H}$, on a small neighborhood $\mathscr{V}_{y_0}\subset X_{H}$ of $y_0$, choose
a unitary trivialization $N|_{\mathscr{V}_{y_0}}\cong \mathscr{V}_{y_0}\times \mathbb{C}^n=\{(y,Z)\,|\,y\in\mathscr{V}_{y_0}, Z=(z_1,\cdots,z_n)\in \mathbb{C}^n\}$
such that for $t\in \mathbb{R}$,
\begin{equation*}
\exp(t\mathcal{H})\cdot\frac{\partial}{\partial z_i}=e^{2\pi\sqrt{-1}\lambda_i t}\frac{\partial}{\partial z_i}.
\end{equation*}
Without loss of generality, we assume that $\lambda_i<0$ for $i\leq \nu$ and $\lambda_i>0$ for $\nu<i\leq n$.
For any $T>0$, $\overrightarrow{k}=(k_1,\cdots, k_n)\in\mathbb{N}^n $, and
$(y,Z)\in \mathscr{V}_{y_0}\times \mathbb{C}^n$, set
\begin{align*}
f_{T,\overrightarrow{k}}(Z)&=\big(\prod_{i=1}^{\nu}z_i^{k_i}\big)\big(\prod_{i=\nu+1}^{n}
\bar{z}_i^{k_i}\big)e^{-T\sum_{i=1}^n\pi|\lambda_i||z_i|^2},
\\
\alpha_{T,\overrightarrow{k}}(y)&=\int_{N_{\R,y}}\rho^2_{\varepsilon}(Z)\prod_{i=1}^n
\left(|z_i|^{2k_i}e^{-2T\pi|\lambda_i||z_i|^2}\right)
\frac{dv_N}{(2\pi)^{\dim N_{\R}}}.
\end{align*}
Computing directly, we have for $s\in{\rm E}^0_{\xi}$ that (compare with \cite[Proposition 9.2]{MR1188532})
\begin{equation}
\begin{split}\label{y9}
p_{T,\xi}s(y,Z)=&\sum_{\overrightarrow{k},\xi_2,\
s.t.\ \sum_{i=1}^{n}k_i|\lambda_i|+\xi_2=\xi}\alpha_{T,\overrightarrow{k}}^{-1}
\rho_{\varepsilon}(Z)f_{T,\overrightarrow{k}}(Z)
\\
&\quad\quad\cdot q_{\xi_2}\int_{N_{\R,y}}\rho_{\varepsilon}(Z')\overline{f_{T,\overrightarrow{k}}(Z')}
s(y,Z')\frac{dv_N(Z')}{(2\pi)^{\dim N_{\R}}}.
\end{split}
\end{equation}

Using (\ref{y9}), we get the following analogue of \cite[Proposition 9.3]{MR1188532}.
\begin{prop}\label{bo6}
There exists $C>0$ such that if $T\geq 1$, $\sigma\in {\rm F}_{\xi}^1$, then
\begin{equation}\label{bo7}
\|I_{T,\xi}\sigma\|_{{\rm \bf{E}}_\xi^1}\leq
C(\|\sigma\|_{{\rm F}_\xi^1}+\sqrt{T}\|\sigma\|_{{\rm F}_\xi^0}).
\end{equation}

There exists $C>0$ such that for any $T\geq 1$, any $s\in
{\rm \bf E}_\xi^1$, then
\begin{equation}\label{bo8}
\|p_{T,\xi}s\|_{{\rm \bf E}_\xi^1}\leq
C(\|s\|_{{\rm \bf E}_\xi^1}+\sqrt{T}\|s\|_{{\rm \bf E}_\xi^0}).
\end{equation}

Given $\gamma>0$, there exists $C'>0$ such that for $T\geq 1$, for
$s\in{{\rm \bf E}_\xi^0}$, then
\begin{equation}\label{bo9}
\|p_{T,\xi}|Z|^\gamma s\|_{{\rm \bf E}_\xi^0}\leq
\frac{C'}{T^{\frac{\gamma}{2}}}\|s\|_{{\rm \bf E}_\xi^0}.
\end{equation}
\end{prop}

\

Since we have the identification of the bundles
\begin{equation*}
\begin{split}
(S(U,L)&\otimes W)|_{\mathscr{V}_{\varepsilon_0}}\simeq
\\
&\left.\widetilde{\pi}^*\left(S(TX_H\oplus V^{\mathbb{R}}_0,L_F)\otimes\Lambda(\overline{N}^*)\otimes
\bigl(\widehat{\otimes}_{v\neq 0}\Lambda V_v\otimes W\bigr)|_{X_H}\right)\right|_{\mathscr{B}_{\varepsilon_0}},
\end{split}
\end{equation*}
we can consider $k^{-1/2}I_{T,\xi}\sigma$ as
an element of ${\rm E}^p_\xi$ for $\sigma\in {\rm F}^p_\xi$. Set
\begin{equation}\label{jtxi}
J_{T,\xi}=k^{-1/2}I_{T,\xi}.
\end{equation}
We denote by $J_{T,\xi,\partial}: {\rm F}_{\xi,\partial}^p\rightarrow {\rm E}_{\xi,\partial}^p$ the restriction of $J_{T,\xi}$ on the boundary.
Let ${\rm E}^p_{T,\xi}=J_{T,\xi}{\rm F}^p_\xi$ (resp. ${\rm E}^p_{T,\xi,\partial}=J_{T,\xi,\partial}{\rm F}^p_{\xi,\partial}$) be the image of
$J_{T,\xi}$ (resp. $J_{T,\xi,\partial}$). Denote the orthogonal complement of ${\rm E}_{T,\xi}^0$ (resp. ${\rm E}_{T,\xi,\partial}^0$) in
${\rm E}_\xi^0$ (resp. ${\rm E}_{\xi,\partial}^0$) by ${\rm E}_{T,\xi}^{0,\bot}$ (resp. ${\rm E}_{T,\xi,\partial}^{0,\bot}$) and let
${\rm E}_{T,\xi}^{p,\bot}={\rm E}_{\xi}^{p}\cap {\rm E}_{T,\xi}^{0,\bot}$ (resp.
${\rm E}_{T,\xi,\partial}^{p,\bot}={\rm E}_{\xi,\partial}^{p}\cap {\rm E}_{T,\xi,\partial}^{0,\bot}$). Let
$\bar{p}_{T,\xi}$ (resp. $\bar{p}_{T,\xi,\partial}$) and $\bar{p}_{T,\xi}^\bot$ (resp. $\bar{p}_{T,\xi,\partial}^\bot$) be the orthogonal
projections from ${\rm E}_\xi^0$ (resp. ${\rm E}_{\xi,\partial}^0$) to ${\rm E}_{T,\xi}^0$ (resp.
${\rm E}_{T,\xi,\partial}^0$) and ${\rm E}_{T,\xi}^{0,\bot}$ (resp. ${\rm E}_{T,\xi,\partial}^{0,\bot}$)
respectively. It is clear that
$\bar{p}_{T,\xi}=k^{-1/2}p_{T,\xi}k^{1/2}$ (resp. $\bar{p}_{T,\xi,\partial}=k^{-1/2}p_{T,\xi,\partial}k^{1/2}$).

For any (possibly unbounded) operator $A$ (resp. $B$) on ${\rm E}^0_\xi$ (resp. ${\rm E}^0_{\xi,\partial}$),
we write
\begin{equation}
A=\begin{pmatrix}
A^{(1)} & A^{(2)}
\\
A^{(3)} & A^{(4)}
\end{pmatrix}
\hspace{20pt}
(\hspace{5pt}
\text{resp.}
\hspace{5pt}
B=\begin{pmatrix}
B^{(1)} & B^{(2)}
\\
B^{(3)} & B^{(4)}
\end{pmatrix}
\hspace{5pt})
\end{equation}
according to the decomposition ${\rm E}^0_\xi={\rm E}^0_{T,\xi}\bigoplus
{\rm E}_{T,\xi}^{0,\bot}$ (resp. ${\rm E}^0_{\xi,\partial}={\rm E}^0_{T,\xi,\partial}\bigoplus
{\rm E}_{T,\xi,\partial}^{0,\bot}$), i.e.,
\begin{align}\label{uy14}
A^{(1)}&=\bar{p}_{T,\xi} A\,\bar{p}_{T,\xi}\ ,\quad
A^{(2)}=\bar{p}_{T,\xi} A\,\bar{p}_{T,\xi}^\bot\ ,
\\
A^{(3)}&=\bar{p}_{T,\xi}^\bot A\,\bar{p}_{T,\xi}\ ,\quad
A^{(4)}=\bar{p}_{T,\xi}^\bot A\,\bar{p}_{T,\xi}^\bot\ .
\notag\\
\notag\\\label{uy15}
(\hspace{5pt}\text{resp.}\hspace{35pt}
B^{(1)}&=\bar{p}_{T,\xi,\partial} B\,\bar{p}_{T,\xi,\partial}\ ,\quad
B^{(2)}=\bar{p}_{T,\xi,\partial} B\,\bar{p}_{T,\xi,\partial}^\bot\ ,
\\
B^{(3)}&=\bar{p}_{T,\xi,\partial}^\bot B\,\bar{p}_{T,\xi,\partial}\ ,\quad
B^{(4)}=\bar{p}_{T,\xi,\partial}^\bot B\,\bar{p}_{T,\xi,\partial}^\bot\ .\notag
\hspace{15pt})
\end{align}

For $T>0$, set
\begin{equation}\label{uy21}
D_{T}=D^X+\sqrt{-1}Tc(H),\quad
B_T=B^X-\sqrt{-1}Tc\left(\frac{\partial}{\partial r}\right)c(H).
\end{equation}
Then $D_T$ is a Dirac type operator with its canonical boundary operator $B_T$
in the sense of Definition \ref{11y1}. Let $D_{T,\xi}$ and $B_{T,\xi}$ be the restrictions of $D_T$ and $B_T$ on ${\rm E}^0_\xi$ and ${\rm E}^0_{\xi,\partial}$,
respectively.

We now introduce a deformation of $D_{T,\xi}$ (resp. $B_{T,\xi}$) according to the decomposition \eqref{uy14} (resp. \eqref{uy15}).
\begin{defn}{\rm (cf. \cite[Definition 3.2]{MR1870658}, \cite[(1.39)]{MR2016198})}\label{aps3}
For any $T>0$, $u\in [0,1]$, set
\begin{equation}
\begin{split}\label{uy22}
D_{T,\xi}(u)&=D_{T,\xi}^{(1)}+D_{T,\xi}^{(4)}+u\bigl(D_{T,\xi}^{(2)}+D_{T,\xi}^{(3)}\,\bigr)\,,
\\
B_{T,\xi}(u)&=B_{T,\xi}^{(1)}+B_{T,\xi}^{(4)}+u\bigl(B_{T,\xi}^{(2)}+B_{T,\xi}^{(3)}\,\bigr)\,.
\end{split}
\end{equation}
\end{defn}

One verifies that $B_{T,\xi}(u)$ is the canonical boundary operator associated to $D_{T,\xi}(u)$ in the sense of
\eqref{liu8}.

\subsection{Various estimates of the operators as $T\rightarrow +\infty$}\label{sec3.4}
We continue the discussion in the previous subsection.
Corresponding to the involution $\tau$ on $S(U,L)$, for $\tau=\tau_s$ (resp. $\tau=\tau_e$), let $D^{X_H}_\xi$ be the restriction of the twisted Spin$^c$ Dirac operator $D^{X_H}\otimes R_{+}(1)$ (resp. $D^{X_H}\otimes R_{-}(1)$)
on ${\rm F}^0_\xi$, and let $B^{X_H}_\xi$ be the restriction of the canonical boundary operator associated to $D^{X_H}\otimes R_{+}(1)$ (resp. $D^{X_H}\otimes R_{-}(1)$) on ${\rm F}^0_{\xi,\partial}$.

With \eqref{y9}, \eqref{uy21} and Propositions \ref{liu33}, \ref{liu45}, \ref{bo6} at our hands, by proceeding
exactly as in \cite[Sections 8 and 9]{MR1188532}, we can show that the following estimates for $B_{T,\xi}^{(i)}$ ($1\leq i\leq 4$) hold.
\begin{prop}{\rm (Compare with \cite[Proposition 3.3]{MR1870658})}\label{lbo3} There exists $\varepsilon >0$ such that
\begin{enumerate}[{\rm (i)}]

\item As $T\longrightarrow +\infty$,
\begin{equation}\label{ly3}
J_{T,\xi,\partial}^{-1} B_{T,\xi}^{(1)}  J_{T,\xi,\partial}=B^{X_H}_\xi
+O\left(\frac{1}{\sqrt{T}}\right),
\end{equation}
where $O(\tfrac{1}{\sqrt{T}})$ denotes a first order differential
operator whose coefficients are dominated by
$\tfrac{C}{\sqrt{T}}\ (C>0)$.

\item
There exist $C_1>0$, $C_2>0$ and $T_0>0$ such that for any $T\geq T_0$, any
$s\in {\rm E}_{T,\xi,\partial}^{1,\bot}$\ ,
$s'\in {\rm E}_{T,\xi,\partial}^1$\ , then
\begin{equation}
\begin{split}
\bigl\|B_{T,\xi}^{(2)}s\bigr\|_{{\rm E}_{\xi,\partial}^0}
&\leq
C_1\Bigl(\frac{1}{\sqrt{T}}\|s\|_{{\rm E}_{\xi,\partial}^1}+\|s\|_{{\rm E}_{\xi,\partial}^0}\Bigr),
\\
\bigl\|B_{T,\xi}^{(3)}s'\bigr\|_{{\rm E}_{\xi,\partial}^0}
&\leq
C_1\Bigl(\frac{1}{\sqrt{T}}\|s'\|_{{\rm E}_{\xi,\partial}^1}+\|s'\|_{{\rm E}_{\xi,\partial}^0}\Bigr),
\end{split}
\end{equation}
and
\begin{equation}
\bigl\|B_{T,\xi}^{(4)}s\bigr\|_{{\rm E}_{\xi,\partial}^0}
\geq C_2\Bigl(\|s\|_{{\rm E}_{\xi,\partial}^1}+\sqrt{T}\,\|s\|_{{\rm E}_{\xi,\partial}^0}\Bigr).
\end{equation}
\end{enumerate}
\end{prop}

\

From here, by proceeding as in \cite[Section 3(c)]{MR1870658}, we can deduce that there exists $T_1>0$
such that for $T\geq T_1$, each $B_{T,\xi}(u)$, for $u\in[0,1]$, is self-adjoint, elliptic and has
discrete eigenvalues with finite multiplicity. Let $P_{T,\xi}(u)$ denote the Atiyah-Patodi-Singer
projection associated to $B_{T,\xi}(u)$.

For any $T\geq T_1$ and $u\in[0,1]$, let
\begin{equation*}
D_{{\rm APS},T,\xi}(u):\bigl\{s\in {{\rm E}^1_\xi}\,\big|\,P_{T,\xi}(u)(s|_{\partial X})=0\bigr\}\longrightarrow {{\rm E}^0_\xi}
\end{equation*}
be the uniquely determined extension of $D_{T,\xi}(u)$.

\begin{prop}{\rm (Compare with \cite[Proposition 3.5]{MR1870658})}\label{aps13}
There exists $T_2>0$ such that for any $u\in [0,1]$ and $T\geq T_2$,
$D_{{\rm APS},T,\xi}(u)$ is a Fredholm operator.
\end{prop}

To prove Proposition \ref{aps13}, we modify the process in \cite[Section 3(d)]{MR1870658}.
For the case where $s$ is supported in
$X\backslash \mathscr{U}_{\epsilon'}$ ($0<\epsilon'<\epsilon$), we
need an analogue of \cite[Lemma 3.7]{MR1870658}. As a matter of fact,
using \eqref{y9}, \eqref{uy21} as well as Propositions \ref{liu33}, \ref{liu45}, \ref{bo6} and proceeding exactly
as in \cite[Sections 8 and 9]{MR1188532},
we deduce the following interior estimates.
\begin{prop}\label{lbo1} There exists $\varepsilon >0$ such that
\begin{enumerate}[{\rm (i)}]
\item As $T\longrightarrow +\infty$,
\begin{equation}\label{ly4}
J_{T,\xi}^{-1}\,D_{T,\xi}^{(1)}J_{T,\xi}=D^{X_H}_\xi+O\left(\frac{1}{\sqrt{T}}\right),
\end{equation}
where $O(\tfrac{1}{\sqrt{T}})$ denotes a first order differential
operator whose coefficients are dominated by
$\tfrac{C}{\sqrt{T}}$ $(C>0)$.

\item There exist $C_1'>0$, $C_2'>0$ and $T_0'>0$ such that for any $T\geq T_0'$, any
$s\in {\rm E}_{T,\xi}^{1,\bot}$\ ,
$s'\in {\rm E}_{T,\xi}^1$\ , then
\begin{equation}
\begin{split}\label{lbo2}
\bigl\|D_{T,\xi}^{(2)}s\bigr\|_{{\rm E}_{\xi}^0}
&\leq
C_1'\left(\frac{\|s\|_{{\rm E}_{\xi}^1}}{\sqrt{T}}+\|s\|_{{\rm E}_{\xi}^0}\right),
\\
\bigl\|D_{T,\xi}^{(3)}s'\bigr\|_{{\rm E}_{\xi}^0}
&\leq
C_1'\left(\frac{\|s'\|_{{\rm E}_{\xi}^1}}{\sqrt{T}}+\|s'\|_{{\rm E}_{\xi}^0}\right),
\end{split}
\end{equation}
and
\begin{equation}
\bigl\|D_{T,\xi}^{(4)}s\bigr\|_{{\rm E}_{\xi}^0}
\geq C_2'\Bigl(\|s\|_{{\rm E}_{\xi}^1}+\sqrt{T}\,\|s\|_{{\rm E}_{\xi}^0}\,\Bigr).
\end{equation}
\end{enumerate}
\end{prop}

With Proposition \ref{lbo1} at our hands, we can complete the proof of Proposition \ref{aps13} in the same way
as in the proof of \cite[Proposition 3.5]{MR1870658} by applying the gluing argument in \cite[pp. 115-117]{MR1188532}.

\subsection{A proof of Theorem \ref{liu28}}\label{sec3.5}
Let $D^{Y_H}_\xi$ be the induced
operator from $B^{X_H}_{\xi}$ through
$\pi_{X_H}$. We first assume that $D^{Y_H}_{\xi}$ is invertible,
then $B^{X_H}_{\xi}$ is invertible. Moreover,
we have the following analogue of \cite[Proposition 3.8]{MR1870658}.

\begin{prop}\label{dz3.8}
If $D^{Y_H}_{\xi}$ is invertible, then
there exists $T_3 >0$ such that for any $T \geq T_3$, $u\in [0,1]$, the boundary operator $B_{T,\xi}(u)$ is invertible.
\end{prop}

By Propositions \ref{aps13} and \ref{dz3.8}, we have a
continuous family of Fredholm operators
$\{D_{{\rm APS},T,\xi}(u)\}_{0 \leq u\leq
1}$ when $T$ is large enough. Furthermore, by Proposition
\ref{dz3.8} and Green's formula, we know that the operators
$D_{{\rm APS},T,\xi}(u)$, ${0 \leq u\leq 1}$, are
self-adjoint. By the homotopy invariance of the index of Fredholm operators, we get
\begin{equation}\label{hopy1}
\tr\left[\tau\big|_{\ker (D_{{\rm APS},T,\xi}(0))}\right]=\tr\left[\tau\big|_{\ker (D_{{\rm APS},T,\xi}(1))}\right].
\end{equation}

\begin{thm}\label{lbo4}{\rm (Compare with \cite[(1.43)]{MR2016198})}
If $D^{Y_H}_{\xi}$ is invertible, then there exists $T_4>0$
such that for any $T \geq T_4$, the following identity holds,
\begin{equation}\label{uy24}
\APS(D_{T,\xi})=\sum_\alpha(-1)^{\sum_{0<v}\dim N_v}\APS\bigl(D_\xi^{X_{H,\alpha}}\bigr)\ .
\end{equation}
\end{thm}

\noindent{\bf Proof}\hspace{8pt} By the definitions of $D_{{\rm APS},T,\xi}(u)$ and $D_{T,\xi}(u)$, we get that
\begin{equation}\label{liu48}
\APS(D_{T,\xi})=\APS(D_{T,\xi}(1))=\tr\left[\tau\big|_{\ker (D_{{\rm APS},T,\xi}(1))}\right].
\end{equation}

Let $P_{T,\xi,1}$ (resp. $P_{T,\xi,4}$) be the Atiyah-Patodi-Singer projection
associated to $B_{T,\xi}^{(1)}$ (resp. $B_{T,\xi}^{(4)}$)
acting on ${\rm E}^0_{T,\xi,\partial}$ (resp. ${\rm E}^{0,\bot}_{T,\xi,\partial}$).
Let
\begin{equation*}
\begin{split}
D_{{\rm APS},T,\xi}^{(1)}&:\bigl\{s\in {{\rm E}^1_{T,\xi}}\,\big|\,P_{T,\xi,1}(s|_{\partial X})=0\bigr\}\longrightarrow {{\rm E}^0_{T,\xi}},
\\
D_{{\rm APS},T,\xi}^{(4)}&:\bigl\{s\in {{\rm E}^{1,\bot}_{T,\xi}}\,\big|\,P_{T,\xi,4}(s|_{\partial X})=0\bigr\}\longrightarrow {{\rm E}^{0,\bot}_{T,\xi}}
\end{split}
\end{equation*}
be the uniquely determined extensions of $D_{T,\xi}^{(1)}$ and $D_{T,\xi}^{(4)}$, respectively.
Using Proposition \ref{lbo3} and proceeding as in the proof of \cite[Proposition 3.5]{MR1870658},
one sees that for $T$ large enough, $D_{{\rm APS},T,\xi}^{(1)}$
and $D_{{\rm APS},T,\xi}^{(4)}$ are both self-adjoint Fredholm operators.
Furthermore,
we  deduce that for $T$ large enough, $\ker (D_{{\rm APS},T,\xi}^{(4)})=0$. Thus we get
\begin{equation}\label{hopy3}
\tr\left[\tau\big|_{\ker (D_{{\rm APS},T,\xi}(0))}\right]=\tr\left[\tau\big|_{\ker (D^{(1)}_{{\rm APS},T,\xi})}\right].
\end{equation}

On the other hand, for $T$ large enough and $u\in [0,1]$, set
\begin{equation}
\begin{split}
D^{X_H}_{T,\xi}(u)&=u\,D^{X_H}_\xi+(1-u)\,J_{T,\xi}^{-1}\,D_{T,\xi}^{(1)}\,J_{T,\xi}\ ,
\\
B^{X_H}_{T,\xi}(u)&=u\,B^{X_H}_\xi+(1-u)\,J_{T,\xi,\partial}^{-1}\,B_{T,\xi}^{(1)}\,J_{T,\xi,\partial}\ .
\end{split}
\end{equation}
From (\ref{ly3}), one can proceed as in \cite[(3.37)-(3.39)]{MR1870658} to see
that when $T$ is large enough, $B^{X_H}_{T,\xi}(u)$ is invertible for every $u\in [0,1]$.

We denote by $P^{X_H}_{T,\xi}(u)$ the Atiyah-Patodi-Singer projection associated to
$B^{X_H}_{T,\xi}(u)$. Using (\ref{ly3}), (\ref{ly4}) and applying the same gluing argument
\cite[pp. 115-117]{MR1188532} as in the proof of
\cite[Proposition 3.5]{MR1870658}, one sees
that when $T$ is large enough and $u\in [0,1]$,
\begin{equation*}
D^{X_H}_{{\rm APS},T,\xi}(u):\bigl\{s\in {{\rm F}^1_\xi}\,\big|\,P^{X_H}_{T,\xi}(u)(s|_{\partial X})=0\bigr\}\longrightarrow {{\rm F}^0_\xi},
\end{equation*}
the uniquely determined extensions of $D^{X_H}_{T,\xi}(u)$,
form a continuous family of self-adjoint Fredholm operators. Thus by the
homotopy invariance of the index of Fredholm operators, one gets
\begin{equation}\label{hopy2}
\tr\left[\tau\big|_{\ker (D^{X_H}_{{\rm APS},T,\xi}(0))}\right]=\tr\left[\tau\big|_{\ker (D^{X_H}_{{\rm APS},T,\xi}(1))}\right]=\APS\bigl(D_\xi^{X_{H}}\bigr).
\end{equation}

From (\ref{liu17}), \eqref{liu18}, (\ref{itxi}) and (\ref{jtxi}), one gets
\begin{equation}\label{ly9}
J^{-1}_{T,\xi}\circ\tau\circ J_{T,\xi}=(-1)^{\sum_{0<v}\dim N_v}\tau,\quad\text{where }\
\tau=\tau_s\text{ or }\tau_e\ .
\end{equation}

From \eqref{hopy1} and \eqref{liu48}-\eqref{ly9}, one finds
\begin{equation}
\APS(D_{T,\xi})=\sum_\alpha(-1)^{\sum_{0<v}\dim N_v}\APS\bigl(D_\xi^{X_{H,\alpha}}\bigr)\ .
\end{equation}

The proof of Theorem \ref{lbo4} is completed.

\

In general, $\dim \ker\bigl(D^{Y_H}_{\xi}\bigr)$ need not be zero.
For any $\xi\in\mathbb{Z}$, choose $a_\xi>0$ be such
that
\begin{equation}\label{uy31}
{\rm Spec}(D^{Y_H}_{\xi})\cap [-2a_\xi,2a_\xi]\subseteq\{0\}.
\end{equation}
To control the eigenvalues of $B_{T,\xi}$ near zero, we use the method in \cite[Section 4(a)]{MR1870658} to perturb the Dirac operators under consideration.

Let $\epsilon>0$ be sufficiently small so that there exists an $S^1$-invariant smooth function
$f: X\longrightarrow \mathbb{R}$
such that $f\equiv 1$ on $\mathscr{U}_{\epsilon/3}$ and  $f\equiv 0$ outside of $\mathscr{U}_{{2\epsilon}/3}$.

Let $D^{X_H}_{\xi,-a_\xi}$  be the Dirac type operator defined by
\begin{equation}
D^{X_H}_{\xi,-a_\xi}=D^{X_H}_{\xi}-a_\xi f c\left(\frac{\partial}{\partial r}\right),
\end{equation}
where for $\tau=\tau_s$ (resp. $\tau_e$), $D^{X_H}_{\xi,-a_\xi}$ is considered as a
differential operator acting on $\Gamma\bigl(X_H,S(TX_H\oplus V_0^{\R},L_F)\otimes R_{+,\xi}\bigr)$
(resp. $\Gamma\bigl(X_H,S(TX_H\oplus V_0^{\R},L_F)\otimes R_{-,\xi}\bigr)$).

By Theorem \ref{thm1}, we get
\begin{equation}\label{uy16}
\APS\bigl(D^{X_{H,\alpha}}_{\xi,-a_\xi}\bigr)-\APS\bigl(D^{X_{H,\alpha}}_{\xi}\bigr)=
-\SF\bigl\{B^{X_{H,\alpha}}_{\xi,+}-a_\xi t\,\big|\,0\leq t\leq1\bigr\}\ .
\end{equation}
By \eqref{uy31}, the right hand side of \eqref{uy16} is equal to zero.

For any $T\in \mathbb{R}$, let $D_{T,-a_\xi}:\Gamma (X,S(U,L)\otimes W) \longrightarrow \Gamma (X,S(U,L)\otimes W)$
be the Dirac type operator defined by
\begin{equation}
D_{T,-a_\xi}=D_T-a_\xi f c\left(\frac{\partial}{\partial r}\right)\ .
\end{equation}
Let $D_{T,\xi,-a_\xi}$ be its restriction to the weight-$\xi$ subspace.

Let $B^{X_H}_{\xi,-a_\xi}$ be the canonical boundary operator of $D^{X_H}_{\xi,-a_\xi}$
in the sense of \eqref{liu8}. Since $D^{Y_H}_\xi-a_\xi$, which is the induced operator
from $B^{X_H}_{\xi,-a_\xi}$ through $\pi_{X_H}$,
is invertible, by the proof of Theorem \ref{lbo4}, we get when $T$ is large enough,
\begin{equation}\label{uy17}
\APS(D_{T,\xi,-a_\xi})=\sum_\alpha(-1)^{\sum_{0<v}\dim N_v}\APS \bigl(D^{X_{H,\alpha}}_{\xi,-a_\xi}\bigr)\ .
\end{equation}

By Theorem \ref{thm1}, we deduce that, for $\xi\in \mathbb{Z}$,
\begin{equation}\label{uy18}
\APS(D_{T,\xi,-a_\xi})\equiv\APS(D_{T,\xi})\mod k\mathbb{Z}.
\end{equation}

From \eqref{uy16}, \eqref{uy17} and \eqref{uy18}, we get
\begin{equation}\label{uy23}
\APS(D_{T,\xi})\equiv\sum_\alpha(-1)^{\sum_{0<v}\dim N_v}\APS \bigl(D^{X_{H,\alpha}}_{\xi}\bigr) \mod k\mathbb{Z}\ .
\end{equation}

On the other hand, by Theorem \ref{thm1}, one knows the {\rm mod} $k$ invariance of $\APS(D_{T,\xi})$ with respect to $T\in \mathbb{R}$,
from which and \eqref{uy24}, \eqref{uy23}, one gets
\begin{equation}
\APS(D,\xi)\equiv\sum_\alpha(-1)^{\sum_{0<v}\dim N_v}\APS \bigl(D^{X_{H,\alpha}}_{\xi}\bigr) \mod k\mathbb{Z}\ .
\end{equation}

By taking $\tau=\tau_s$ (resp. $\tau_e$), we get the first equation of \eqref{liu29} (resp. \eqref{liu30}).
To get the second equation of \eqref{liu29} (resp. \eqref{liu30}), we only need to apply the first equation of
\eqref{liu29} (resp. \eqref{liu30}) to the case where the circle action on $X$ defined by the inverse of the
original circle action on $X$.

The proof of Theorem \ref{liu28} is completed.

\section{Rigidity and vanishing theorems on $\mathbb{Z}/k$ Spin$^c$ manifolds}\label{sec4}

In this section, combining the $S^1$-equivariant index theorem we have established in Section \ref{sec2}
with the methods of \cite{MR1870666}, we prove the rigidity and vanishing theorems for
$\mathbb{Z}/k$ Spin$^c$ manifolds, which generalize \cite[Theorems A and B]{MR1396769}. As
will be pointed out in Remark \ref{11y2}, when applied to $\mathbb{Z}/k$ spin manifolds, our results
provide a resolution to a conjecture of Devote \cite{MR1405063}.
Both the statement of the main results and their proof are inspired by
the corresponding results as well as their proof for closed manifolds in \cite{MR1870666, MR2016198}.
As explained in Section 2.1, when we regard the considered $\mathbb{Z}/k$ manifold as a quotient space
which has the homotopy type of a CW complex, by using splitting
principle \cite[Chapter 17]{MR1249482}, we can apply  the topological arguments
in \cite{MR1870666, MR2016198} in our $\mathbb{Z}/k$ context with little modification. Thus
we will only indicate the main steps of the proof of our results.

This section is organized as follows. In Section \ref{sec4.1}, we state our main results, the rigidity
and vanishing theorems for $\mathbb{Z}/k$ Spin$^c$ manifolds. In Section \ref{sec4.2}, we present two recursive
formulas which will be used to prove our main results stated in Section \ref{sec4.1}.
In Section \ref{sec4.3}, we prove the rigidity and vanishing theorems for
$\mathbb{Z}/k$ Spin$^c$ manifolds.

\subsection{Rigidity and vanishing theorems}\label{sec4.1}

Let $X$ be a $2l$-dimensional $\mathbb{Z}/k$ manifold, which admits a nontrivial $\mathbb{Z}/k$ circle action.
We assume that $TX$ has a $\mathbb{Z}/k$ $S^1$-equivariant Spin$^c$ structure. Let $V$ be an even dimensional
$\mathbb{Z}/k$ real vector bundle over $X$. We assume that $V$ has a $\mathbb{Z}/k$ $S^1$-equivariant spin structure.
Let $W$ be a $\mathbb{Z}/k$ $S^1$-equivariant complex vector bundle of rank $r$ over $X$.
Let $K_W=\det(W)$ be the determinant line bundle of $W$, which is obviously
a $\mathbb{Z}/k$ complex line bundle.

Let $K_X$ be the $\mathbb{Z}/k$
complex line bundle over $X$ induced by the Spin$^c$ structure of $TX$.
Let $S(TX, K_X)$ be the complex spinor bundle of $(TX, K_X)$ as in Section \ref{sec2.1}.
Let $S(V)=S^+(V)\oplus S^-(V)$ be the spinor bundle of $V$.

Let $K(X)$ be the $K$-group of $\mathbb{Z}/k$ complex vector bundles over $X$ (cf. \cite[pp. 285]{MR1144425}).
We define the following elements in $K(X)[[q^{1/2}]]$ (cf. \cite[(2.1)]{MR1870666})
\begin{equation}\label{eg3}
\begin{split}
R_1(V)&=\Bigl(S^+(V)+S^-(V)\Bigr)\otimes\bigotimes_{n=1}^{\infty}\Lambda_{q^n}(V)\ ,
\\
R_2(V)&=\Bigl(S^+(V)-S^-(V)\Bigr)\otimes\bigotimes_{n=1}^{\infty}\Lambda_{-q^n}(V)\ ,
\\
R_3(V)&=\bigotimes_{n=1}^{\infty}\Lambda_{-q^{n-1/2}}(V)\ ,
\\
R_4(V)&=\bigotimes_{n=1}^{\infty}\Lambda_{q^{n-1/2}}(V)\ .
\end{split}
\end{equation}
For $N\in \mathbb{N}$, let $y=e^{2\pi i/N}\in \mathbb{C}$ be an $N$th root of unity. Set
\begin{equation}\label{eg4}
Q_{y}(W)=\bigotimes_{n=0}^{\infty}\Lambda_{-y^{-1}\cdot q^n}(\overline{W})
\otimes\bigotimes_{n=1}^{\infty}\Lambda_{-y\cdot q^n}(W)\in K(X)[[q]]\ .
\end{equation}
Then there exist $Q_{\ell}(W)\in K(X)[[q]]$, $0\leq\ell< N$ such that
\begin{equation}\label{quan}
Q_{y}(W)=\sum_{\ell=0}^{N-1}y^{\ell}Q_{\ell}(W).
\end{equation}

Let $H^*_{S^1}(X,\mathbb{Z})=H^*(X\times_{S^1}ES^1,\mathbb{Z})$ denote the $S^1$-equivariant cohomology group
of $X$, where $ES^1$ is the universal $S^1$-principal bundle over the classifying space $BS^1$ of $S^1$. So
$H^*_{S^1}(X,\mathbb{Z})$ is a module over $H^*(BS^1,\mathbb{Z})$ induced by the projection
$\overline{\pi}:X\times_{S^1}ES^1\rightarrow BS^1$.
Let $p_1(\cdot)_{S^1}$ and $\omega_2(\cdot)_{S^1}$ denote the first $S^1$-equivariant pontrjagin class and
the second $S^1$-equivariant Stiefel-Whitney class, respectively. As $V\times_{S^1}ES^1$ is spin over $X\times_{S^1}ES^1$, one knows that $\frac{1}{2}p_1(V)_{S^1}$ is well defined in $H^*_{S^1}(X,\mathbb{Z})$ (cf. \cite[pp. 456-457]{MR998662}). Recall that
\begin{equation}\label{eg5}
H^*(BS^1,\mathbb{Z})=\mathbb{Z}[[u]]
\end{equation}
with $u$ a generator of degree 2.

In the following, we denote by $D^X\otimes R$ the twisted Spin$^c$ Dirac operator acting on $S(TX,K_X)\otimes R$ (cf. Definition \ref{diracdef}).
Furthermore, for $m\in\frac{1}{2}\mathbb{Z}$, $h\in \mathbb{Z}$ and
$R(q)=\sum_{m\in\frac{1}{2}\mathbb{Z}}q^mR_m\in K_{S^1}(X)[[q^{1/2}]]$, we
will also denote $\APS(D^X\otimes R_m,h)$ (cf. \eqref{ly5}) by $\APS(D^X\otimes R(q),m,h)$.

Now we can state the main results of this paper as follows, which generalize \cite[Theorems A and B]{MR1396769}
to the case of $\mathbb{Z}/k$ Spin$^c$ manifolds.
\begin{thm}\label{maintheorem}
Assume that $\omega_2(W)_{S^1}=\omega_2(TX)_{S^1}$, $\frac{1}{2}p_1(V+W-TX)_{S^1}=e\cdot\overline{\pi}^*u^2$ $(e\in\mathbb{Z})$ in $H^*_{S^1}(X,\mathbb{Z})$, and $c_1(W)=0 \mod (N)$. For
$0\leq\ell< N$, $i=1$, $2$, $3$, $4$, consider the $S^1$-equivariant twisted Spin$^c$ Dirac operators
\begin{equation}
D^X\otimes(K_W\otimes K_X^{-1})^{1/2}\otimes\bigotimes_{n=1}^{\infty}{\rm Sym}_{q^n}(TX)\otimes R_i(V)\otimes Q_{\ell}(W)\ .
\end{equation}
\begin{enumerate}[{\rm (i)}]
\item If $e=0$, then these operators are rigid in $\mathbb{Z}/k$ category.

\item If $e<0$, then they have vanishing properties in $\mathbb{Z}/k$ category.
\end{enumerate}
\end{thm}

\begin{rem}{\rm (Compare with \cite[Remark 2.1]{MR1870666})}
As $\omega_2(W)_{S^1}=\omega_2(TX)_{S^1}$, $c_1(K_W\otimes K_X^{-1})_{S^1}=0 \mod(2)$. We note that in our case,
$X\times_{S^1}ES^1$ has the homotopy type of a CW complex {\rm \cite{MR0290354}}. By {\rm \cite[Corollary 1.2]{MR0461538}},
the circle action on $X$ can be lifted to $(K_W\otimes K_X^{-1})^{1/2}$ and is compatible with the circle action on $K_W\otimes K_X^{-1}$.
\end{rem}

\begin{rem}\label{11y2}
If X is a $\mathbb{Z}/k$ spin manifold, by taking $V=TX$, $W=0$ and $i=3$ in Theorem \ref{maintheorem}, we resolve a conjecture of {\rm \cite{MR1405063}}.
\end{rem}

Actually, as in \cite{MR1870666}, our proof of Theorem \ref{maintheorem} works under the following slightly weaker hypothesis.
Let us first explain some notations.

For each $n>1$, consider $\mathbb{Z}_n\subset S^1$, the cyclic subgroup of order $n$. We have the $\mathbb{Z}_n$-equivariant
cohomology of $X$ defined by $H^*_{\mathbb{Z}_n}(X,\mathbb{Z})=H^*(X\times_{\mathbb{Z}_n}ES^1,\mathbb{Z})$,
and there is a natural \textquoteleft\textquoteleft forgetful\textquoteright\textquoteright\ map $\alpha(S^1,\mathbb{Z}_n):X\times_{\mathbb{Z}_n}ES^1\rightarrow X\times_{S^1}
ES^1$ which induces a pullback $\alpha(S^1,\mathbb{Z}_n)^*:H^*_{S^1}(X,\mathbb{Z})\rightarrow
H^*_{\mathbb{Z}_n}(X,\mathbb{Z})$. We
denote by $\alpha(S^1,1)$ the arrow which forgets the $S^1$-action.
Thus $\alpha(S^1,1)^*:H^*_{S^1}(X,\mathbb{Z})\rightarrow H^*(X,\mathbb{Z})$
is induced by the inclusion of $X$ into $X\times_{S^1}ES^1$ as a fiber over $BS^1$.

Finally, note that if $\mathbb{Z}_n$ acts trivially on a space $M$, then there is a new arrow $t^*:H^*(M,\mathbb{Z})
\rightarrow H^*_{\mathbb{Z}_n}(M,\mathbb{Z})$ induced by the projection
$t:M\times_{\mathbb{Z}_n}ES^1=M\times B\mathbb{Z}_n\rightarrow M$.

Let $\mathbb{Z}_{\infty}=S^1$. For each $1<n\leq+\infty$, let $i:X(n)\rightarrow X$ be the
inclusion of the fixed point set of $\mathbb{Z}_n\subset S^1$ in $X$, and so $i$ induces
$i_{S^1}:X(n)\times_{S^1}ES^1\rightarrow X\times_{S^1}ES^1$.

In the rest of this paper, we use the same assumption as in \cite[(2.4)]{MR1870666}.
Suppose that there exists some
integer $e\in\mathbb{Z}$ such that for $1<n\leq+\infty$,
\begin{equation}\label{hypothesis}
\begin{split}
\alpha(S^1,\mathbb{Z}_n)^*\circ i_{S^1}^*\Bigl(\,\frac{1}{2}p_1(V+W-TX)_{S^1}-e\cdot\overline{\pi}^*u^2\Bigr)&
\\
=t^*\circ\alpha(S^1,1)^*\circ i_{S^1}^*\Bigl(\,\frac{1}{2}p_1(V+W-TX)_{S^1}\Bigr)\ .&
\end{split}
\end{equation}

Remark that the relation (\ref{hypothesis}) clearly follows from the hypothesis of Theorem \ref{maintheorem} by
pulling back and forgetting. Thus it is a weaker hypothesis.

Let $G_y$ be the multiplicative group generated by $y$. Following Witten \cite{MR970288},  we consider the
action of $y_0\in G_y$ on $W$ (resp. $\overline{W}$) by multiplication by $y_0$ (resp. $y_0^{-1}$) on $W$ (resp. $\overline{W}$).
Set \begin{equation}\label{eg6}
Q(W)=\bigotimes_{n=0}^{\infty}\Lambda_{-q^n}(\overline{W})
\otimes\bigotimes_{n=1}^{\infty}\Lambda_{-q^n}(W)\in K(X)[[q]]\ .
\end{equation}
Then the actions of $G_y$ on $W$ and $\overline{W}$ naturally induce the action of $G_y$ on $Q(W)$.
Clearly, $y\cdot Q(W)=Q_y(W)$.
By (\ref{quan}), we know that for $0\leq \ell<N$,
\begin{equation}
y_0\cdot Q_\ell(W)=y_0^\ell Q_\ell(W),\quad\ \text{where }\ y_0\in G_y.
\end{equation}

In what follows, for $m\in\frac{1}{2}\mathbb{Z}$, $0\leq\ell<N$, $h\in \mathbb{Z}$ and
$R(q)\in K_{S^1}(X)[[q^{1/2}]]$, we will denote $\APS(D^X\otimes R(q)\otimes Q_\ell(W),m,h)$
by $\APS(D^X\otimes R(q)\otimes Q(W),m,\ell,h)$.

We can now state a slightly more general version of Theorem \ref{maintheorem}.
\begin{thm}\label{theorem}
Under the hypothesis {\rm (\ref{hypothesis})}, consider the $S^1\times G_y$-equivariant twisted Spin$^c$ Dirac operators
\begin{equation}
D^X\otimes(K_W\otimes K_X^{-1})^{1/2}\otimes\bigotimes_{n=1}^{\infty}{\rm Sym}_{q^n}(TX)\otimes R_i(V)\otimes Q(W)\ .
\end{equation}
\begin{enumerate}[{\rm (i)}]
\item If $e=0$, for $m\in \frac{1}{2}\mathbb{Z}$, $h\in\mathbb{Z}$, $h\neq 0$, $0\leq\ell<N$, one has
\begin{equation}\label{index}
\begin{split}
\APS\Bigl(D^{X}&\otimes(K_W\otimes K_X^{-1})^{1/2}\otimes\bigotimes_{n=1}^{\infty}{\rm Sym}_{q^n}(TX)\Bigr.
\\
\Bigl.&\otimes R_i(V)\otimes Q(W),m,\ell,h\Bigr)\equiv 0 \mod k\mathbb{Z}\ .
\end{split}
\end{equation}

\item If $e<0$, for $m\in \frac{1}{2}\mathbb{Z}$, $h\in\mathbb{Z}$, $0\leq\ell<N$, one has
\begin{equation}
\begin{split}
\APS\Bigl(D^{X}&\otimes(K_W\otimes K_X^{-1})^{1/2}\otimes\bigotimes_{n=1}^{\infty}{\rm Sym}_{q^n}(TX)\Bigr.
\\
\Bigl.&\otimes R_i(V)\otimes Q(W),m,\ell,h\Bigr)\equiv 0 \mod k\mathbb{Z}\ .
\end{split}
\end{equation}
In particular, one has
\begin{equation}
\begin{split}
\APS\Bigl(D^{X}&\otimes(K_W\otimes K_X^{-1})^{1/2}\otimes\bigotimes_{n=1}^{\infty}{\rm Sym}_{q^n}(TX)\Bigr.
\\
\Bigl.&\otimes R_i(V)\otimes Q(W),m,\ell\Bigr)\equiv 0 \mod k\mathbb{Z}\ .
\end{split}
\end{equation}
\end{enumerate}
\end{thm}

\

The rest of this paper is devoted to a proof of Theorem \ref{theorem}.

\subsection{Several intermediate results}\label{sec4.2}

Recall that $X_H=\{X_{H,\alpha}\}$ be the fixed point set of the circle action. As in \cite[pp. 940]{MR1870666},
we may and we will assume that
\begin{equation}\label{eg8}
\begin{split}
&TX|_{X_H}=TX_H\oplus\bigoplus_{v>0}N_v,\\
&TX|_{X_H}\otimes_{\mathbb{R}}\mathbb{C}=TX_H\otimes_{\mathbb{R}}\mathbb{C}
\oplus\bigoplus_{v>0}\,\bigl(N_v\oplus\overline{N}_v\bigr),
\end{split}
\end{equation}
where $N_v$ is the complex vector bundles on which $S^1$ acts by sending $g$ to $g^v$. We also assume that
\begin{equation}\label{eg9}
V|_{X_H}=V_0^\mathbb{R}\oplus\bigoplus_{v>0}V_v\ ,\qquad
W|_{X_H}=\bigoplus_{v}W_v\ ,
\end{equation}
where $V_v$, $W_v$ are complex vector bundles on which $S^1$ acts by sending $g$ to $g^v$, and $V_0^{\mathbb{R}}$
is a real vector bundle on which $S^1$ acts as identity.

By \eqref{eg8}, as in \eqref{liu17} or \eqref{liu18}, there is a natural $\mathbb{Z}/k$ isomorphism
between the $\mathbb{Z}_2$-graded $C(TX)$-Clifford modules over $X_H$,
\begin{equation}\label{eg10}
S(TX,K_X)|_{X_H}\simeq S\Bigl(TX_H,K_X\otimes_{v>0}(\det N_v)^{-1}\Bigr)
\,\widehat{\otimes}\,\widehat{\bigotimes}_{v>0}\Lambda N_v\ .
\end{equation}

For a $\mathbb{Z}/k$ complex vector bundle $R$ over $X_H$, let $D^{X_H}\otimes R$,
$D^{X_{H,\alpha}}\otimes R$ be the twisted Spin$^c$ Dirac operators on
$S(TX_{H},K_X\otimes_{v>0}(\det N_v)^{-1})\otimes R$ over $X_H$, $X_{H,\alpha}$, respectively (cf. Definition \ref{diracdef}).

For $i=1$, $2$, $3$, $4$, we set
\begin{equation}
R^i=(K_W\otimes K_X^{-1})^{1/2}\otimes R_i(V)\otimes Q(W) .
\end{equation}

Then by Theorem \ref{liu28}, we can express the global Atiyah-Patodi-Singer index
via the Atiyah-Patodi-Singer indices on the fixed point set up to $k\mathbb{Z}$.

\begin{prop}\label{eg1}{\rm (Compare with \cite[Proposition 2.1]{MR1870666})}
For $m\in \frac{1}{2}\mathbb{Z}$, $h\in \mathbb{Z}$, $1\leq i\leq 4$, $0\leq\ell<N$, we have
\begin{equation}\label{eg2}
\begin{split}
&\APS\,\Bigl(D^X\otimes\otimes_{n=1}^{\infty}{\rm Sym}_{q^n}(TX)\otimes R^i,m,\ell,h\Bigr)
\\
&\equiv\sum_{\alpha}(-1)^{\sum_{v>0}\dim N_v}
\APS\,\Bigl(D^{X_{H,\alpha}}\otimes_{n=1}^{\infty}{\rm Sym}_{q^n}(TX)\otimes R^i\Bigr.
\\
\Bigl.&\hspace{90pt}\otimes{\rm Sym}\,(\oplus_{v>0}N_v)\otimes_{v>0}\det N_v,m,\ell,h\Bigr)\quad {\rm mod}\,k\mathbb{Z}\ .
\end{split}
\end{equation}
\end{prop}

To simplify the notations, we use the same convention as in \cite[pp. 945]{MR1870666}. For $n_0\in \mathbb{N}^*$, we define a number operator
$P$ on $K_{S^1}(X)[[q^{\frac{1}{n_0}}]]$ in the following way: if $R(q)=\oplus_{n\in\frac{1}{n_0}\mathbb{Z}}
R_nq^n\in K_{S^1}(X)[[q^{\frac{1}{n_0}}]]$, then $P$ acts on $R(q)$ by multiplication by $n$ on $R_n$. From
now on, we simply denote ${\rm Sym}_{q^n}(TX)$, $\Lambda_{q^n}(V)$ and $\Lambda_{q^n}(W)$
by ${\rm Sym}(TX_n)$, $\Lambda(V_n)$ and $\Lambda(W_n)$, respectively. In this way, $P$ acts on $TX_n$, $V_n$ and $W_n$
by multiplication by $n$, and the actions of $P$ on ${\rm Sym}(TX_n)$, $\Lambda(V_n)$ and $\Lambda(W_n)$ are naturally
induced by the corresponding actions of $P$ on $TX_n$, $V_n$ and $W_n$. So the eigenspace of $P=n$ is just
given by the coefficient of $q^n$ of the corresponding element $R(q)$. For $R(q)=\oplus_{n\in\frac{1}{n_0}\mathbb{Z}}
R_nq^n\in K_{S^1}(X)[[q^{\frac{1}{n_0}}]]$, we will also denote
$\APS\bigl(D^{X}\otimes R_m,h\bigr)$ by $\APS\bigl(D^{X}\otimes R(q),m,h\bigr)$.

Recall that $H$ is the Killing vector field on $X$ corresponding to $\mathcal{H}$, the canonical basis of ${\rm Lie}(S^1)$. If $E$ is a $\mathbb{Z}/k$ $S^1$-equivariant vector bundle over $X$, let $\mathscr{L}_H$
denote the corresponding Lie derivative along $H$ acting on $\Gamma(X_H,E|_{X_H})$. The weight of the circle action on
$\Gamma(X_H,E|_{X_H})$ is given by the action
\begin{equation*}
\emph{\textbf{J}}_H=\frac{1}{2\pi\sqrt{-1}}\mathscr{L}_H.
\end{equation*}
Recall that the $\mathbb{Z}_2$-grading on $S(TX,K_X)\otimes_{n=1}^{\infty}{\rm Sym}(TX_n)$
is induced by
the $\mathbb{Z}_2$-grading on $S(TX,K_X)$. Write
\begin{equation}\label{eg12}
\begin{split}
&\mathcal{F}^0(X)=\bigotimes_{n=1}^{\infty}{\rm Sym}_{q^n}(TX)\otimes{\rm Sym}(\oplus_{v>0}N_v)\otimes_{v>0}\det N_v\ ,
\\
&F_V^1=S(V)\otimes\bigotimes_{n=1}^{\infty}\Lambda(V_n)\ ,\qquad
F_V^2=\bigotimes_{n\in\mathbb{N}+\frac{1}{2}}\Lambda(V_n)\ ,
\\
&Q^1(W)=\bigotimes_{n=0}^\infty \Lambda(\overline{W}_n)\otimes \bigotimes_{n=1}^\infty\Lambda(W_n)\ .
\end{split}
\end{equation}
There are two natural $\mathbb{Z}_2$-gradings on $F_V^1$, $F_V^2$ (resp. $Q^1(W)$). The first grading is induced by the $\mathbb{Z}_2$-grading
of $S(V)$ and the forms of homogeneous degrees in $\otimes_{n=1}^{\infty}\Lambda(V_n)$,
$\otimes_{n\in\mathbb{N}+\frac{1}{2}}\Lambda(V_n)$ (resp. $Q^1(W)$). We define $\tau_e|_{F_V^{i\pm}}=\pm 1$ ($i=1$, $2$)
(resp. $\tau_1|_{Q^1(W)^{\pm}}=\pm1$) to be the involution
defined by this $\mathbb{Z}_2$-grading. The second grading is the one for which $F_V^i$ ($i=1$, $2$) are purely even,
i.e., $F_V^{i+}=F_V^i$. We denote by $\tau_s={\rm id}$ the involution defined by this $\mathbb{Z}_2$-grading.
Then the coefficient of $q^n$ ($n\in\frac{1}{2}\mathbb{Z}$) in (\ref{eg3}) of $R_1(V)$ or $R_2(V)$
(resp. $R_3(V)$, $R_{4}(V)$ or $Q(W)$) is exactly the $\mathbb{Z}_2$-graded $\mathbb{Z}/k$ vector subbundle of $(F_V^1,\tau_s)$
or $(F_V^1,\tau_e)$ (resp. $(F_V^2,\tau_e)$, $(F_V^2,\tau_s)$ or $(Q^1(W),\tau_1)$), on which $P$ acts by multiplication by $n$.

Furthermore, we denote by $\tau_e$ (resp. $\tau_s$) the $\mathbb{Z}_2$-grading on $S(TX,K_X)
\otimes\otimes_{n=1}^{\infty}{\rm Sym}(TX_n)\otimes F_V^i$ ($i=1$, $2$) induced by
the above $\mathbb{Z}_2$-gradings. We will denote by $\tau_{e1}$ (resp. $\tau_{s1}$)
the $\mathbb{Z}_2$-grading on $S(TX,K_X)\otimes\otimes_{n=1}^{\infty}{\rm Sym}(TX_n)\otimes F_V^i\otimes Q^1(W)$ ($i=1$, $2$) defined by
\begin{equation}\label{eg13}
\tau_{e 1}=\tau_e\widehat{\otimes}\tau_1\ ,\quad \tau_{s 1}=\tau_s\widehat{\otimes}\tau_1\ .
\end{equation}

Let $h^{V_v}$ be the Hermitian metric on $V_v$ induced by the metric $h^V$ on $V$. In the following, we identity
$\Lambda V_v$ with $\Lambda \overline{V}^*_v$ by using the Hermitian metric $h^{V_v}$ on $V_v$. By (\ref{eg9}),
as in (\ref{eg10}), there is a natural $\mathbb{Z}/k$ isomorphism between the
$\mathbb{Z}_2$-graded $C(V)$-Clifford modules over $X_H$,
\begin{equation}\label{eg14}
S(V)|_{X_H}\simeq S\Bigl(V_0^\mathbb{R},\otimes_{v>0}(\det V_v)^{-1}\Bigr)\otimes\widehat{\bigotimes}_{v>0}\Lambda V_v\ .
\end{equation}

Let $V_0=V_0^{\mathbb{R}}\otimes_{\mathbb{R}}\mathbb{C}$. By using the above notations,
we rewrite (\ref{eg12}) on the fixed point set $X_H$,
\begin{equation}\label{eg15}
\begin{split}
&\mathcal{F}^0(X)=\bigotimes_{n=1}^{\infty}{\rm Sym}(TX_H)\otimes\bigotimes_{n=1}^{\infty}{\rm Sym}
\Bigl(\oplus_{v>0}\bigl(N_{v,n}\oplus\overline{N}_{v,n}\bigr)\Bigr)
\\
&\hspace{70pt}\otimes{\rm Sym}(\oplus_{v>0}N_v)\otimes_{v>0}\det N_v\ ,
\\
&F_V^1=\bigotimes_{n=1}^{\infty}\Lambda\Bigl(V_{0,n}\oplus\oplus_{v>0}(V_{v,n}\oplus\overline{V}_{v,n})\Bigr)
\\
&\hspace{70pt}\otimes S\Bigl(V_0^\mathbb{R},\otimes_{v>0}(\det V_v)^{-1}\Bigr)\otimes_{v>0}\Lambda V_{v,0}\ ,
\\
&F_V^2=\bigotimes_{n\in\mathbb{N}+\frac{1}{2}}\Lambda\Bigl(V_{0,n}\oplus\oplus_{v>0}(V_{v,n}\oplus\overline{V}_{v,n})\Bigr)\ ,
\\
&Q^1(W)=\bigotimes_{n=0}^{\infty}\Lambda(\oplus_{v}\overline{W}_{v,n})
\otimes\bigotimes_{n=1}^{\infty}\Lambda(\oplus_{v}W_{v,n})\ .
\end{split}
\end{equation}

We introduce the same shift operators as in \cite[Section 3.2]{MR1870666}, which follow \cite{MR998662} in spirit.
For $p\in \mathbb{N}$, we set
\begin{equation}\label{eg16}
\begin{split}
&r_*:N_{v,n}\rightarrow N_{v,n+pv}\ ,\quad r_*:\overline{N}_{v,n}\rightarrow \overline{N}_{v,n-pv}\ ,
\\
&r_*:V_{v,n}\rightarrow V_{v,n+pv}\ ,\quad \ r_*:\overline{V}_{v,n}\rightarrow \overline{V}_{v,n-pv}\ ,
\\
&r_*:W_{v,n}\rightarrow W_{v,n+pv}\ ,\quad r_*:\overline{W}_{v,n}\rightarrow \overline{W}_{v,n-pv}\ .
\end{split}
\end{equation}

Furthermore, for $p\in \mathbb{N}$, we introduce the following elements in
$K_{S^1}(X_H)[[q]]$ (cf. \cite[(3.6)]{MR1870666}),
\begin{equation}\label{eg21}
\begin{split}
&\mathcal{F}_p(X)=\bigotimes_{n=1}^{\infty}{\rm Sym}(TX_H)\otimes\bigotimes_{v>0}\Bigl(\,\bigotimes_{n=1}^{\infty}
{\rm Sym}\,(N_{v,n})\bigotimes_{n> pv}{\rm Sym}\,(\overline{N}_{v,n})\,\Bigr)\ ,
\\
&\mathcal{F}_p'(X)=\bigotimes_{v>0}\bigotimes_{0\leq n\leq pv}\Bigl({\rm Sym}\,(N_{v,-n})\otimes\det (N_v)\Bigr)\ ,
\\
&\mathcal{F}^{-p}(X)=\mathcal{F}_p(X)\otimes\mathcal{F}_p'(X)\ .
\end{split}
\end{equation}
Note that when $p=0$, $\mathcal{F}^{-p}(X)$ is exactly the $\mathcal{F}^{0}(X)$ in (\ref{eg15}).
The $\mathbb{Z}_2$-grading on $S(TX_H, K_X\otimes_{v>0}(\det N_v)^{-1})\otimes \mathcal{F}^{-p}(X)\,$
is induced by
the $\mathbb{Z}_2$-grading on $S(TX_H, K_X\otimes_{v>0}(\det N_v)^{-1})$.

As in \cite[(2.9)]{MR1870666}, we write
\begin{equation}\label{eg17}
\begin{split}
L(N)&=\otimes_{v>0}(\det N_v)^v\ ,\quad L(V)=\otimes_{v>0}(\det V_v)^v\ ,
\\
L(W)&=\otimes_{v\neq0}(\det W_v)^v\ ,\quad
L=L(N)^{-1}\otimes L(V)\otimes L(W)\ .
\end{split}
\end{equation}

Using the similar $\mathbb{Z}/k$ $S^1$-equivariant isomorphism
of complex vector bundles as in \cite[(3.14)]{MR2016198}
and the similar $\mathbb{Z}/k$ $G_y\times S^1$-equivariant isomorphism of complex vector bundles as in
\cite[(3.15) and (3.16)]{MR1870666}, by direct calculation, we deduce the following proposition.
\begin{prop}\label{eg22}{\rm (cf. \cite[Proposition 3.1]{MR1870666})}
For $p\in \mathbb{Z}$, $p>0$, $i=1$, $2$, there are natural $\mathbb{Z}/k$ isomorphisms of vector bundles over $X_H$,
\begin{equation}\label{eg23}
r_*(\mathcal{F}^{-p}(X))\simeq\mathcal{F}^{0}(X)\otimes L(N)^{p}\ ,\quad
r_*(F_V^i)\simeq F_V^i\otimes L(V)^p\ .
\end{equation}
For any $p\in \mathbb{Z}$, $p>0$, there is a natural $\mathbb{Z}/k$ $G_y\times S^1$-equivariant isomorphism of vector bundles over $X_H$,
\begin{equation}\label{eg24}
r_*(Q^1(W))\simeq Q^1(W)\otimes L(W)^{-p}.
\end{equation}
\end{prop}

On $X_H$, as in \cite[(2.8)]{MR1870666}, we write
\begin{equation}
\begin{split}
e(N)&=\sum_{v>0}v^2\dim N_v\ ,\quad d'(N)=\sum_{v>0}v\dim N_v\ ,
\\
e(V)&=\sum_{v>0}v^2\dim V_v\ ,\quad d'(V)=\sum_{v>0}v\dim V_v\ ,
\\
e(W)&=\sum_{v}v^2\dim W_v\ ,\quad d'(W)=\sum_{v}v\dim W_v\ .
\end{split}
\end{equation}
Then $e(N)$, $e(V)$, $e(W)$, $d'(N)$, $d'(V)$ and $d'(W)$ are locally constant functions on $X_H$.

Take $\mathbb{Z}_{\infty}=S^1$ in the hypothesis (\ref{hypothesis}).  By using splitting principle \cite[Chapter 17]{MR1249482}, we get the same identities as in \cite[(2.11)]{MR1870666},
\begin{equation}\label{uy50}
c_1(L)=0\ ,\quad e(V)+e(W)-e(N)=2e\ .
\end{equation}
As indicated in Section \ref{sec2.1}, \eqref{uy50} means $L$ is a trivial complex line bundle over
each component $X_{H,\alpha}$ of $X_H$, and $S^1$ acts on $L$ by sending $g$ to $g^{2e}$, and $G_y$ acts on $L$
by sending $y$ to $y^{d'(W)}$.

The following proposition is deduced from Proposition \ref{eg22}.
\begin{prop}\label{eg25}{\rm (cf. \cite[Proposition 3.2]{MR1870666})}
For $p\in \mathbb{Z}$, $p>0$, $i=1$, $2$, the $\mathbb{Z}/k$
$G_y$-equivariant isomorphism of vector bundles over $X_H$
induced by {\rm (\ref{eg23})}, {\rm (\ref{eg24})},
\begin{equation}\label{eg26}
\begin{split}
r_*\,:&\ S(TX_H,K_X\otimes_{v>0}(\det N_v)^{-1})\otimes(K_W\otimes K^{-1}_X)^{1/2}
\\
&\hspace{50pt}\otimes\mathcal{F}^{-p}(X)
\otimes F^i_V\otimes Q^1(W)
\\
&\longrightarrow\ S(TX_H,K_X\otimes_{v>0}(\det N_v)^{-1})\otimes(K_W\otimes K^{-1}_X)^{1/2}
\\
&\hspace{50pt}\otimes\mathcal{F}^{0}(X)
\otimes F^i_V\otimes Q^1(W)\otimes L^{-p}\ ,
\end{split}
\end{equation}
verifies the following identities
\begin{equation}\label{eg27}
\begin{split}
r_*^{-1}\cdot \textbf{J}_H\cdot r_*&=\textbf{J}_H\ ,
\\
r_*^{-1}\cdot P\cdot r_*&=P+p\textbf{J}_H+p^2e-\frac{1}{2}p^2e(N)-\frac{p}{2}d'(N)\ .
\end{split}
\end{equation}
For the $\mathbb{Z}_2$-gradings, we have
\begin{equation}\label{eg28}
r_*^{-1}\tau_{e}r_*=\tau_{e}\ ,
\quad r_*^{-1}\tau_{s}r_*=\tau_{s}\ ,
\quad r_*^{-1}\tau_{1}r_*=(-1)^{pd'(W)}\tau_{1}\ .
\end{equation}
\end{prop}

Then we get the following recursive formula.
\begin{thm}\label{eg29}
{\rm (Compare with \cite[Theorem 2.4]{MR1870666})} For each $\alpha$,
$m\in \tfrac{1}{2}\mathbb{Z}$, $1\leq \ell<N$, $h$, $p\in \mathbb{Z}$, $p>0$, the following identity holds,
\begin{align}\label{eg30}
&\APS\Bigl(D^{X_{H,\alpha}}\otimes {\cal F}^{-p}(X)\otimes R^i,m+\frac{1}{2}p^2e(N)+\frac{p}{2}d'(N),\ell,h\Bigr)
\\
&=(-1)^{pd'(W)}\APS\Bigl(D^{X_{H,\alpha}}\otimes\mathcal{F}^{0}(X)\otimes R^i\otimes L^{-p},m+ph+p^2e,\ell,h\Bigr)\ .
\notag
\end{align}
\end{thm}

Now we state another recursive formula whose proof will be presented in Section 5.

\begin{thm}\label{eg31}
{\rm (Compare with \cite[Theorem 2.3]{MR1870666})} For $1\leq i\leq 4$, $m\in \tfrac{1}{2}\mathbb{Z}$, $1\leq \ell<N$,
$h$, $p\in \mathbb{Z}$, $p>0$, we have the following identity,
\begin{equation}\label{eg32}
\begin{split}
&\sum_{\alpha}(-1)^{\sum_{v>0}\dim N_v}\APS\Bigl(D^{X_{H,\alpha}}\otimes\mathcal{F}^0(X)\otimes R^i,m,\ell,h\Bigr)
\\
&\equiv\sum_{\alpha}(-1)^{pd'(N)+\sum_{v>0}\dim N_v}\APS\Bigl(D^{X_{H,\alpha}}\otimes\mathcal{F}^{-p}(X)
\otimes R^i,\bigr.
\\
&\Bigl.\hspace{75pt}
m+\tfrac{1}{2}p^2e(N)+\tfrac{p}{2}d'(N),\ell,h\Bigr) \mod k\mathbb{Z}\ .
\end{split}
\end{equation}
\end{thm}

\subsection{A proof of Theorem \ref{theorem}}\label{sec4.3}

Recall we assume in Theorem \ref{maintheorem} that $c_1(W)= 0$ mod $(N)$. Then by \cite[Section 8]{MR981372} and
\cite[Lemma 2.1]{MR1870666}, $d'(W)$ mod$(N)$ is constant on each connected component $X_{H,\alpha}$ of $X_H$. So we can
extend $L$ to a trivial complex line bundle over $X$, and we extend the $S^1$-action on it by sending $g$
on the canonical section 1 of $L$ to $g^{2e}\cdot1$, and $G_y$ acts on $L$ by sending $y$ to $y^{d'(W)}$.

As $\tfrac{1}{2}p_1(TX-W)_{S^1}\in H^*_{S^1}(X,\mathbb{Z})$ is well defined, one has the same identity as in \cite[(2.27)]{MR1870666},
\begin{equation}\label{eg33}
d'(N)+d'(W)\equiv0 \mod (2).
\end{equation}

From Proposition \ref{eg1}, Theorems \ref{eg29}, \ref{eg31} and (\ref{eg33}), for $1\leq i\leq 4$, $m\in\tfrac{1}{2}\mathbb{Z}$,
$1\leq \ell<N$, $h$, $p\in\mathbb{Z}$, $p>0$, we get the following identity
(compare with \cite[(2.28)]{MR1870666}),
\begin{align}\label{eg34}
&\APS\bigl(D^{X}\otimes\bigotimes_{n=1}^{\infty}{\rm Sym}_{q^n}(TX)\otimes R^i,m,\ell,h\bigr)
\\
&\equiv\APS\bigl(D^{X}\otimes\bigotimes_{n=1}^{\infty}{\rm Sym}_{q^n}(TX)\otimes R^i
\otimes L^{-p},m',\ell,h\bigr) \mod k\mathbb{Z}\ ,\notag
\end{align}
with
\begin{equation}\label{eg35}
m'=m+ph+p^2e.
\end{equation}

By (\ref{eg3}), (\ref{eg4}), if $m<0$ or $m'<0$, then either side of (\ref{eg34}) is identically zero,
which completes the proof of  Theorem \ref{theorem}. In fact,
\begin{enumerate}[{\rm (i)}]
\item Assume that $e=0$. Let $h\in \mathbb{Z}$, $m_0\in\tfrac{1}{2}\mathbb{Z}$, $h\neq 0$ be fixed. If $h>0$, we take $m'=m_0$, then for $p$ large enough,
we get $m<0$ in (\ref{eg34}). If $h<0$, we take $m=m_0$, then for $p$ large enough,
we get $m'<0$ in (\ref{eg34}).

\item Assume that $e<0$. For $h\in\mathbb{Z}$, $m_0\in\tfrac{1}{2}\mathbb{Z}$, we take $m=m_0$, then for $p$ large enough,
we get $m'<0$ in (\ref{eg34}).
\end{enumerate}

The proof of  Theorem \ref{theorem} is completed.

\section{A proof of Theorem \ref{eg31}}\label{sec5}

In this section, following \cite[Section 4]{MR1870666}, we present a proof of Theorem \ref{eg31}.

This section is organized as follows. In Section \ref{sec5.1}, we first introduce the same refined shift operators as in \cite[Section 4.2]{MR1870666}.
In Section \ref{sec5.2}, we construct the twisted Spin$^c$ Dirac operator on $X(n_j)$, the fixed point
set of the naturally induced $\mathbb{Z}_{n_j}$-action on $X$. In Section \ref{sec5.3}, by applying the
$S^1$-equivariant index theorem we have established in Section \ref{sec2}, we prove Theorem \ref{eg31}.

\subsection{The refined shift operators}\label{sec5.1}

We first introduce a partition of $[0,1]$ as in \cite[pp. 942--943]{MR1870666}. Set
$J=\bigl\{v\in \mathbb{N}\,\big| \text{ there exists } \alpha \text{ such that } N_v\neq0 \text{ on } X_{H,\alpha}\bigr\}$ and
\begin{equation}\label{eg39}
\Phi=\bigl\{\beta\in (0,1]\,\big| \text{ there exists } v\in J \text{ such that } \beta v\in\mathbb{Z}\bigr\}\ .
\end{equation}
We order the elements in $\Phi$ so that $\Phi=\bigl\{\beta_i\,\big|\,1\leq i\leq J_0, J_0\in\mathbb{N}\text{ and }
\beta_i<\beta_{i+1}\bigr\}$. Then for any integer $1\leq i\leq J_0$, there exist $p_i$, $n_i\in\mathbb{N}$,
$0<p_i\leq n_i$, with $(p_i,n_i)=1$ such that
\begin{equation}\label{eg40}
\beta_i=p_i/n_i\ .
\end{equation}
Clearly, $\beta_{J_0}=1$. We also set $p_0=0$ and $\beta_0=0$.

For $0\leq j\leq J_0$, $p\in\mathbb{N}^*$, we write
\begin{equation}\label{eg41}
\begin{split}
I_j^p&=\Bigl\{(v,n)\in\mathbb{N}\times\mathbb{N}\,\Big|\,v\in J,\ (p-1)v<n\leq pv,
\ \frac{n}{v}=p+1+\frac{p_j}{n_j}\Bigr\}\ ,
\\
\overline{I}_j^p&=\Bigl\{(v,n)\in\mathbb{N}\times\mathbb{N}\,\Big|\,v\in J,\ (p-1)v<n\leq pv,\
\frac{n}{v}>p+1+\frac{p_j}{n_j}\Bigr\}\ .
\end{split}
\end{equation}
Clearly, $I_0^p=\emptyset$, the empty set. We define $\mathcal{F}_{p,j}(X)$ as in \cite[(2.21)]{MR1870666},
which are analogous with (\ref{eg21}). More specifically, we set
\begin{align}
&\mathcal{F}_{p,j}(X)=\bigotimes_{n=1}^{\infty}{\rm Sym}\,(TX_H)\otimes\bigotimes_{v>0}\Biggl(\,\bigotimes_{n=1}^{\infty}
{\rm Sym}\,(N_{v,n})\otimes\bigotimes_{n> (p-1)v+\tfrac{p_j}{n_j}v}{\rm Sym}\,(\overline{N}_{v,n})\,\Biggr)
\notag\\
&\hspace{65pt}\bigotimes_{v>0}\bigotimes_{0\leq n\leq (p-1)v+\bigl[\tfrac{p_j}{n_j}v\bigr]}
\Bigl({\rm Sym}\,(N_{v,-n})\otimes\det N_v\Bigr)\\
&=\mathcal{F}_{p}(X)\otimes\mathcal{F}_{p-1}'(X)\otimes\bigotimes_{(v,n)\in\cup_{i=0}^jI_i^p}
\Bigl({\rm Sym}\,(N_{v,-n})\otimes\det N_v\Bigr)\bigotimes_{(v,n)\in\overline{I}_j^p}{\rm Sym}\,(\overline{N}_{v,n})\ ,\notag
\end{align}
where we use the notation that for $s\in \mathbb{R}$, $[s]$ denotes the greatest integer which is less than or equal to $s$. Then
\begin{equation}
\mathcal{F}_{p,0}(X)=\mathcal{F}^{-p+1}(X)\ ,\quad
\mathcal{F}_{p,J_0}(X)=\mathcal{F}^{-p}(X)\ .
\end{equation}

From the construction of $\beta_i$, we know that for $v\in J$, there is no integer in
$\bigl(\tfrac{p_{j-1}}{n_{j-1}}v,\tfrac{p_j}{n_j}v\bigr)$. Furthermore,
\begin{equation}
\begin{split}
\Bigl[\frac{p_{j-1}}{n_{j-1}}v\Bigr]&=\Bigl[\frac{p_{j}}{n_{j}}v\Bigr]-1\ \text{  if  }\ v\equiv0\mod(n_j)\ ,
\\
\Bigl[\frac{p_{j-1}}{n_{j-1}}v\Bigr]&=\Bigl[\frac{p_{j}}{n_{j}}v\Bigr]\ \text{ if }\ v\not \equiv0\mod(n_j)\ .
\end{split}
\end{equation}

We use the same shift operators $r_{j*}$, $1\leq j\leq J_0$ as in \cite[(4.21)]{MR1870666}, which
refine the shift operator $r_{*}$ defined in (\ref{eg16}). For $p\in\mathbb{N}^*$, set
\begin{equation}\label{eg43}
\begin{split}
&r_{j*}:N_{v,n}\rightarrow N_{v,n+(p-1)v+p_jv/n_j}\ ,\quad
r_{j*}:\overline{N}_{v,n}\rightarrow \overline{N}_{v,n-(p-1)v-p_jv/n_j}\ ,
\\
&r_{j*}:V_{v,n}\rightarrow V_{v,n+(p-1)v+p_jv/n_j}\ ,
\quad \ r_{j*}:\overline{V}_{v,n}\rightarrow \overline{V}_{v,n-(p-1)v-p_jv/n_j}\ ,
\\
&r_{j*}:W_{v,n}\rightarrow W_{v,n+(p-1)v+p_jv/n_j}\ ,
\quad r_{j*}:\overline{W}_{v,n}\rightarrow \overline{W}_{v,n-(p-1)v-p_jv/n_j}\ .
\end{split}
\end{equation}

For $1\leq j\leq J_0$, we define $\mathcal{F}(\beta_j)$, $F_V^1(\beta_j)$, $F_V^2(\beta_j)$ and $Q_W(\beta_j)$ as in \cite[(4.13)]{MR1870666}.
\begin{equation*}
\begin{split}
\mathcal{F}(\beta_j)&=\bigotimes_{0<n\in\mathbb{Z}}{\rm Sym}\,(TX_{H,n})\otimes\bigotimes_{v>0,v\equiv0,\tfrac{n_j}{2}
\,{\rm mod}(n_j) }\bigotimes_{0<n\in\mathbb{Z}+\tfrac{p_j}{n_j}v}{\rm Sym}\,(N_{v,n}\oplus\overline{N}_{v,n})
\\
&\otimes\bigotimes_{0<v'<n_j/2}{\rm Sym}\,\Biggl(\ \bigoplus_{v\equiv v',-v'\ {\rm mod}(n_j)}\Bigl(\bigoplus_
{0<n\in\mathbb{Z}+\tfrac{p_j}{n_j}v}N_{v,n}\bigoplus_{0<n\in\mathbb{Z}-\tfrac{p_j}{n_j}v}
\overline{N}_{v,n}\Bigr)\ \Biggr)\ ,
\end{split}
\end{equation*}
\begin{equation*}
\begin{split}
F_V^1(\beta_j)&=\Lambda\left(\ \bigoplus_{0<n\in\mathbb{Z}}V_{0,n}\bigoplus_{v>0,v\equiv0,\tfrac{n_j}{2}
\,{\rm mod}(n_j)}\Bigl(\bigoplus_
{0<n\in\mathbb{Z}+\tfrac{p_j}{n_j}v}V_{v,n}\bigoplus_{0<n\in\mathbb{Z}-\tfrac{p_j}{n_j}v}
\overline{V}_{v,n}\Bigr)\right.
\\
&\left.\bigoplus_{0<v'<n_j/2}\Biggl(\ \bigoplus_{v\equiv v',-v'\,{\rm mod}(n_j)}\Bigl(\bigoplus_
{0<n\in\mathbb{Z}+\tfrac{p_j}{n_j}v}V_{v,n}\bigoplus_{0<n\in\mathbb{Z}-\tfrac{p_j}{n_j}v}
\overline{V}_{v,n}\Bigr)\Biggr)\right)\ ,
\end{split}
\end{equation*}
\begin{align}\label{uy4}
F_V^2(\beta_j)\!&=\!\Lambda\negthickspace\left(\bigoplus_{\negthickspace 0<n\in\mathbb{Z}+\tfrac{1}{2}}V_{0,n}\bigoplus_{\negthickspace v>0,\,v\equiv0,\tfrac{n_j}{2}
\,{\rm mod}(n_j)}\Bigl(\negthickspace\bigoplus_{\negthickspace 0<n\in\mathbb{Z}+\tfrac{p_j}{n_j}v+\tfrac{1}{2}}V_{v,n}
\negthickspace\bigoplus_{\negthickspace 0<n\in\mathbb{Z}-\tfrac{p_j}{n_j}v+\tfrac{1}{2}}
\overline{V}_{v,n}\!\Bigr)\right.
\notag\\
&\left.\bigoplus_{0<v'<n_j/2}\Biggl(\ \bigoplus_{v\equiv v',-v'\,{\rm mod}(n_j)}\Bigl(\bigoplus_
{0<n\in\mathbb{Z}+\tfrac{p_j}{n_j}v+\tfrac{1}{2}}V_{v,n}\bigoplus_{0<n\in\mathbb{Z}-\tfrac{p_j}{n_j}v+\tfrac{1}{2}}
\overline{V}_{v,n}\Bigr)\Biggr)\right)\ ,
\notag\\
Q_W&(\beta_j)=\Lambda\Biggl(\ \bigoplus_v\ \Bigl(\bigoplus_
{0<n\in\mathbb{Z}+\tfrac{p_j}{n_j}v}W_{v,n}\bigoplus_{0\leq n\in\mathbb{Z}-\tfrac{p_j}{n_j}v}
\overline{W}_{v,n}\ \Bigr)\Biggr)\ .
\end{align}

Using the definition of $r_{j*}$ and computing directly, we get an analogue of Proposition \ref{eg22} as follows.
\begin{prop}\label{eg44}{\rm (cf. \cite[Proposition 4.1]{MR1870666})}
There are natural $\mathbb{Z}/k$ isomorphisms of vector bundles over $X_H$,
\begin{equation*}
\begin{split}
r_{j*}(\mathcal{F}_{p,j-1}(X))&\simeq\mathcal{F}(\beta_j)\otimes\bigotimes_{v>0,\ v\equiv0\ {\rm mod}\,(n_j)}
{\rm Sym}\,(\overline{N}_{v,0})
\\
&\otimes\bigotimes_{v>0}\ (\det N_v)^{\bigl[\tfrac{p_j}{n_j}v\bigr]+(p-1)v+1}
\otimes\bigotimes_{v>0,\ v\equiv0\ {\rm mod}\,(n_j)}(\det N_v)^{-1}\ ,
\\
r_{j*}(\mathcal{F}_{p,j}(X))&\simeq\ \mathcal{F}(\beta_j)\otimes\bigotimes_{v>0,\ v\equiv0\ {\rm mod}\,(n_j)}
{\rm Sym}(N_{v,0})
\\
&\hspace{110pt}
\otimes\bigotimes_{v>0}\ (\det N_v)^{\bigl[\tfrac{p_j}{n_j}v\bigr]+(p-1)v+1}\ ,
\end{split}
\end{equation*}
\begin{equation*}
\begin{split}
r_{j*}(F_V^1)& \simeq S\Bigl(V_0^{\mathbb{R}}\,,\,\otimes_{v>0}(\det V_v)^{-1}\Bigr)\otimes F_V^1(\beta_j)
\\
&\hspace{40pt}\otimes\bigotimes_{v>0,\ v\equiv0\ {\rm mod}\,(n_j)}\Lambda(V_{v,0})\otimes
\bigotimes_{v>0}(\det \overline{V}_v)^{\bigl[\tfrac{p_j}{n_j}v\bigr]+(p-1)v}\ ,
\\
r_{j*}(F_V^2)& \simeq F_V^2(\beta_j)\otimes\bigotimes_{v>0,\ v\equiv\tfrac{n_j}{2}\ {\rm mod}\,(n_j)}
\Lambda(V_{v,0})\otimes\bigotimes_{v>0}(\det \overline{V}_v)^{\bigl[\tfrac{p_j}{n_j}v+\tfrac{1}{2}\bigr]+(p-1)v}\ .
\end{split}
\end{equation*}

There is a natural $\mathbb{Z}/k$ $G_y\times S^1$-equivariant isomorphism of vector bundles over $X_H$,
\begin{equation*}\label{eg46}
\begin{split}
r_{j*}(Q^1(W))&\simeq Q_W(\beta_j)\otimes\bigotimes_{v>0}(\det \overline{W}_v)^{\bigl[\tfrac{p_j}{n_j}v\bigr]+(p-1)v+1}
\\
&\hspace{15pt}
\otimes\bigotimes_{v>0,\ v\equiv0\ {\rm mod}\,(n_j)}(\det\overline{W}_v)^{-1}\otimes
\bigotimes_{v<0}(\det W_v)^{\bigl[-\tfrac{p_j}{n_j}v\bigr]-(p-1)v}\ .
\end{split}
\end{equation*}
\end{prop}

\subsection{The Spin$^c$ Dirac operators on $X(n_j)$}\label{sec5.2}

Recall that there is a nontrivial $\mathbb{Z}/k$ circle action on $X$ which can be lifted to
the $\mathbb{Z}/k$ circle actions on $V$ and $W$.

For $n\in \mathbb{N}^*$, let $\mathbb{Z}_n\subset S^1$
denote the cyclic subgroup of order $n$. Let $X(n_j)$ be the fixed point set of the induced $\mathbb{Z}_{n_j}$
action on $X$. Let $N(n_j)\rightarrow X(n_j)$ be the normal bundle to $X(n_j)$ in $X$. As in \cite[pp. 151]{MR954493}
(see also \cite[Section 4.1]{MR1870666}, \cite[Section 4.1]{MR2016198} or \cite{MR998662}),
we see that $N(n_j)$ and $V$ can be decomposed, as $\mathbb{Z}/k$ real vector bundles over $X(n_j)$, into
\begin{equation}\label{eg48}
\begin{split}
N(n_j)&=\bigoplus_{0<v<n_j/2}N(n_j)_v\oplus N(n_j)_{n_j/2}^{\mathbb{R}}\ ,
\\
V|_{X(n_j)}&= V(n_j)_{0}^{\mathbb{R}}\oplus\bigoplus_{0<v<n_j/2}V(n_j)_v\oplus V(n_j)_{n_j/2}^{\mathbb{R}}\ ,
\end{split}
\end{equation}
where $V(n_j)_{0}^{\mathbb{R}}$ is the $\mathbb{Z}/k$ real vector bundle on which $\mathbb{Z}_{n_j}$ acts by identity,
and $N(n_j)_{n_j/2}^{\mathbb{R}}$ (resp. $V(n_j)_{n_j/2}^{\mathbb{R}}$) is defined to be zero if $n_j$ is odd.
Moreover, for $0<v<n_j/2$, $N(n_j)_v$ (resp. $V(n_j)_v$) admits
unique $\mathbb{Z}/k$ complex structure such that $N(n_j)_v$ (resp. $V(n_j)_v$) becomes a $\mathbb{Z}/k$
complex vector bundle on which $g\in \mathbb{Z}_{n_j}$ acts by $g^v$. We also denote by
$V(n_j)_0$, $V(n_j)_{n_j/2}$ and $N(n_j)_{n_j/2}$ the corresponding complexification of
$V(n_j)_{0}^{\mathbb{R}}$, $V(n_j)_{n_j/2}^{\mathbb{R}}$ and $N(n_j)_{n_j/2}^{\mathbb{R}}$.

Similarly, we also have the following $\mathbb{Z}_{n_j}$-equivariant decomposition of $W$, as
$\mathbb{Z}/k$ complex vector bundles over $X(n_j)$,
\begin{equation}\label{eg49}
W=\bigoplus_{0\leq v<n_j}W(n_j)_v\ ,
\end{equation}
where for $0\leq v<n_j$, $g\in \mathbb{Z}_{n_j}$ acts on $W(n_j)_v$ by sending $g$ to $g^v$.

By \cite[Lemma 4.1]{MR1870666} (see also \cite[Lemmas 9.4 and 10.1]{MR954493} or \cite[Lemma 5.1]{MR998662}),
we know that the $\mathbb{Z}/k$ vector bundles $TX(n_j)$ and $V(n_j)_{0}^{\mathbb{R}}$ are orientable and even dimensional. Thus $N(n_j)$ is orientable
over $X(n_j)$. By (\ref{eg48}), $V(n_j)_{n_j/2}^{\mathbb{R}}$ and $N(n_j)_{n_j/2}^{\mathbb{R}}$ are also
orientable and even dimensional. In what follows, we fix the orientations of $N(n_j)_{n_j/2}^{\mathbb{R}}$
and $V(n_j)_{n_j/2}^{\mathbb{R}}$. Then $TX(n_j)$ and $V(n_j)_0^{\mathbb{R}}$ are naturally oriented by
(\ref{eg48}) and the orientations of $TX$, $V$, $N(n_j)_{n_j/2}^{\mathbb{R}}$ and
$V(n_j)_{n_j/2}^{\mathbb{R}}$.

By (\ref{eg8}), (\ref{eg9}), (\ref{eg48}) and (\ref{eg49}), upon restriction to $X_H$,
we get the following identifications of $\mathbb{Z}/k$ complex vector bundles (cf. \cite[(4.9) and (4.12)]{MR1870666}),
for $0<v\leq n_j/2$,
\begin{equation}
\begin{split}\label{uy2}
N(n_j)_v&=\bigoplus_{v'>0,\,v'\equiv v\,{\rm mod}(n_j)}N_{v'}\ \oplus
\bigoplus_{v'>0,\,v'\equiv -v\,{\rm mod}(n_j)}\overline{N}_{v'}\ \ ,
\\
V(n_j)_v&=\bigoplus_{v'>0,\,v'\equiv v\,{\rm mod}(n_j)}V_{v'}\ \oplus
\bigoplus_{v'>0,\,v'\equiv -v\,{\rm mod}(n_j)}\overline{V}_{v'}\ \ ,
\end{split}
\end{equation}
for $0\leq v<n_j$,
\begin{equation}\label{uy3}
W(n_j)_v=\bigoplus_{v'>0,\,v'\equiv v\,{\rm mod}(n_j)}W_{v'}\ \ .
\end{equation}
Also we get the following identifications of $\mathbb{Z}/k$ real vector bundles over $X_H$
(cf. \cite[(4.11)]{MR1870666}),
\begin{equation*}
\begin{split}
TX(n_j)&=TX_H\oplus\bigoplus_{v>0,\,v\equiv0\,{\rm mod}\,(n_j)}N_v\ ,\quad
N(n_j)_{n_j/2}^{\mathbb{R}}=\bigoplus_
{v>0,\,v\equiv\tfrac{n_j}{2}\,{\rm mod}\,(n_j)}N_v\ \ ,
\\
V(n_j)_0^{\mathbb{R}}&=V_0^{\mathbb{R}}\oplus\bigoplus_{v>0,\,v\equiv0\,{\rm mod}(n_j)}V_v\ ,\quad
V(n_j)_{n_j/2}^{\mathbb{R}}=\bigoplus_{v>0,\,v\equiv\tfrac{n_j}{2}\,{\rm mod}(n_j)}V_v\ \ .
\end{split}
\end{equation*}
Moreover, we have the identifications of $\mathbb{Z}/k$ complex vector bundles over $X_H$ as follows,
\begin{equation}\label{eg51}
\begin{split}
TX(n_j)\otimes_{\mathbb{R}}\mathbb{C}&=TX_H\otimes_{\mathbb{R}}\mathbb{C}\oplus
\bigoplus_{v>0,\,v\equiv0\,{\rm mod}(n_j)}(N_v\oplus\overline{N}_v)\ ,
\\
V(n_j)_0&=V_0^{\mathbb{R}}\otimes_{\mathbb{R}}\mathbb{C}\oplus
\bigoplus_{v>0,\,v\equiv0\,{\rm mod}(n_j)}(V_v\oplus\overline{V}_v)\ .
\end{split}
\end{equation}

As $(p_j,n_j)=1$, we know that, for $v\in\mathbb{Z}$, $p_jv/n_j\in\mathbb{Z}$ if and only if $v/n_j\in\mathbb{Z}$.
Also, $p_jv/n_j\in\mathbb{Z}+\tfrac{1}{2}$ if and only if $v/n_j\in\mathbb{Z}+\tfrac{1}{2}$. Remark if
$v\equiv -v'\ {\rm mod}(n_j)$, then $\{n\,|\,0<n\in\mathbb{Z}+\tfrac{p_j}{n_j}v\}=
\{n\,|\,0<n\in\mathbb{Z}-\tfrac{p_j}{n_j}v'\}$. Using the identifications
(\ref{uy2}), (\ref{uy3}) and (\ref{eg51}), we can rewrite $\mathcal{F}(\beta_j)$, $F_V^1(\beta_j)$,
$F_V^2(\beta_j)$ and $Q_W(\beta_j)$ defined in (\ref{uy4}) as follows (cf. \cite[(4.7)]{MR1870666}),
\begin{equation}
\begin{split}\label{uy8}
\mathcal{F}(\beta_j)&=\bigotimes_{0<n\in\mathbb{Z}}{\rm Sym}\,\bigl(TX(n_j)_n\bigr)\otimes\bigotimes_{0<v<n_j/2}{\rm Sym}\,
\Bigl(\bigoplus_{0<n\in\mathbb{Z}+\tfrac{p_j}{n_j}v}N(n_j)_{v,n}\Bigr.
\\
&\hspace{30pt}
\Bigl.\oplus\bigoplus_{0<n\in\mathbb{Z}-\tfrac{p_j}{n_j}v}
\overline{N(n_j)}_{v,n}\Bigr)\oplus\bigoplus_{0<n\in\mathbb{Z}+\tfrac{1}{2}}{\rm Sym}\,\bigl(N(n_j)_{n_j/2,n}\bigr)\ ,
\end{split}
\end{equation}
\begin{equation}
\begin{split}\label{uy9}
F_V^1(\beta_j)&=\Lambda\Bigl(\ \bigoplus_{0<n\in\mathbb{Z}}V(n_j)_{0,n}\oplus\bigoplus_{0<v<n_j/2}
\Bigl(\ \bigoplus_{0<n\in\mathbb{Z}+\tfrac{p_j}{n_j}v}V(n_j)_{v,n}\Bigr.\Bigr.
\\
&\hspace{30pt}
\Bigl.\Bigl.\oplus\bigoplus_{0<n\in\mathbb{Z}-\tfrac{p_j}{n_j}v}
\overline{V(n_j)}_{v,n}\Bigr)\bigoplus_{0<n\in\mathbb{Z}+\tfrac{1}{2}}V(n_j)_{n_j/2,n}\Bigr)\ ,
\end{split}
\end{equation}
\begin{equation}
\begin{split}\label{uy10}
F_V^2(\beta_j)&=\Lambda\Bigl(\ \bigoplus_{0<n\in\mathbb{Z}}V(n_j)_{n_j/2,n}\oplus\bigoplus_{0<v<n_j/2}
\Bigl(\ \bigoplus_{0<n\in\mathbb{Z}+\tfrac{p_j}{n_j}v+\tfrac{1}{2}}V(n_j)_{v,n}\Bigr.\Bigr.
\\
&\hspace{30pt}
\Bigl.\Bigl.\oplus\bigoplus_{0<n\in\mathbb{Z}-\tfrac{p_j}{n_j}v+\tfrac{1}{2}}\overline{V(n_j)}_{v,n}\Bigr)
\bigoplus_{0<n\in\mathbb{Z}+\tfrac{1}{2}}V(n_j)_{0,n}\Bigr)\ ,
\end{split}
\end{equation}
\begin{equation}\label{eg52}
Q_W(\beta_j)=\Lambda\Bigl(\ \bigoplus_{0\leq v<n_j}\Bigl(\bigoplus_
{0<n\in\mathbb{Z}+\tfrac{p_j}{n_j}v}W(n_j)_{v,n}\oplus\bigoplus_{0\leq n\in\mathbb{Z}-\tfrac{p_j}{n_j}v}
\overline{W(n_j)}_{v,n}\Bigr)\Bigr)\ .
\end{equation}
Thus $\mathcal{F}(\beta_j)$, $F_V^1(\beta_j)$, $F_V^2(\beta_j)$ and $Q_W(\beta_j)$ can be extended to
$\mathbb{Z}/k$ vector bundles over $X(n_j)$.

We now define the Spin$^c$ Dirac operators on $X(n_j)$ following \cite[Section 4.1]{MR1870666}.

Consider the hypothesis in (\ref{hypothesis}).
By splitting principle \cite[Chapter 17]{MR1249482} and computing as in \cite[Lemmas 11.3 and 11.4]{MR954493},
we get
\begin{equation}\label{eg53}
\begin{split}
&\Biggl(\ \sum_{0<v<\tfrac{n_j}{2}}v\cdot c_1\Bigl(V(n_j)_v+W(n_j)_v-W(n_j)_{n_j-v}-N(n_j)_v\Bigr)\Biggr.
\\
&\quad
\Biggl.+r(n_j)\cdot\frac{n_j}{2}\cdot\omega_2\Bigl(W(N_j)_{n_j/2}+V(n_j)_{n_j/2}-N(n_j)_{n_j/2}\Bigr)\Biggr)\cdot u_{n_j}=0\ ,
\end{split}
\end{equation}
where $r(n_j)=\frac{1}{2}(1+(-1)^{n_j})$, and $u_{n_j}\in H^2(B\mathbb{Z}_{n_j},\mathbb{Z})\simeq \mathbb{Z}_{n_j}$
is the generator of $H^*(B\mathbb{Z}_{n_j},\mathbb{Z})\simeq \mathbb{Z}[u_{n_j}]/(n_j\cdot u_{n_j})$\ . Then by (\ref{eg53}), we know that
\begin{equation*}
\begin{split}
&\sum_{0<v<\tfrac{n_j}{2}}v\cdot c_1\Bigl(V(n_j)_v+W(n_j)_v-W(n_j)_{n_j-v}-N(n_j)_v\Bigr)
\\
&\qquad
+r(n_j)\cdot\frac{n_j}{2}\cdot\omega_2\Bigl(W(n_j)_{n_j/2}+V(n_j)_{n_j/2}-N(n_j)_{n_j/2}\Bigr)
\end{split}
\end{equation*}
is divided by $n_j$. Therefore, we have

\begin{lem}\label{eg54}
{\rm (cf. \cite[Lemma 4.2]{MR1870666})} Assume that {\rm (\ref{hypothesis})} holds. Let
\begin{equation}\label{eg55}
\begin{split}
L(n_j)&=\bigotimes_{0<v<n_j/2}\Bigl(\det (N(n_j)_v)\otimes\det(\overline{V(n_j)_v}))\Bigr.
\\&\hspace{60pt}\Bigl.
\otimes\det(\overline{W(n_j)_v})
\otimes\det(W(n_j)_{n_j-v})\Bigr)^{(r(n_j)+1)v}
\end{split}
\end{equation}
be the complex line bundle over $X(n_j)$. Then we have
\begin{enumerate}[{\rm (i)}]
\item $L(n_j)$ has an $n_j^{\rm th}$ root over $X(n_j)$.

\item Let $U_1=TX(n_j)\oplus V(n_j)_0^{\mathbb{R}}$, $U_2=TX(n_j)\oplus V(n_j)_{n_j/2}^{\mathbb{R}}$. Let
\begin{equation*}\label{eg56}
\begin{split}
L_1&=K_X\otimes\bigotimes_{0<v<n_j/2}\Bigl(\det (N(n_j)_v)\otimes\det(\overline{V(n_j)_v})\Bigr)
\\
&\hspace{90pt}
\otimes\det\Bigl(W(n_j)_{n_j/2}\Bigr)
\otimes L(n_j)^{r(n_j)/n_j}\ ,
\\
L_2&=K_X\otimes\bigotimes_{0<v<n_j/2}\Bigl(\det (N(n_j)_v)\Bigr)
\otimes\det\Bigl(W(n_j)_{n_j/2}\Bigr)\otimes L(n_j)^{r(n_j)/n_j}.
\end{split}
\end{equation*}
Then
$U_1$ (resp. $U_2$) has a $\mathbb{Z}/k$ Spin$^c$ structure defined by $L_1$ (resp. $L_2$).
\end{enumerate}
\end{lem}

Remark that in order to define an $S^1$ (resp. $G_y$) action on $L(n_j)^{r(n_j)/n_j}$, we must replace the $S^1$
(resp.  $G_y$) action by its $n_j$-fold action. Here by abusing notation,
we still say an $S^1$ (resp. $G_y$) action without causing any confusion.

In what follows, by $D^{X(n_j)}$ we mean the $S^1$-equivariant Spin$^c$ Dirac operator on $S(U_1,L_1)$ or $S(U_2,L_2)$
over $X(n_j)$ (cf. Definition \ref{diracdef}).

Corresponding to (\ref{liu15}), by (\ref{uy2}),  we denote by
\begin{align}
S(U_1,L_1)'&=S\Bigl(TX_H\oplus V_0^{\mathbb{R}},L_1\otimes\bigotimes_{v>0,\,v\equiv0\,{\rm mod}\,(n_j)}
(\det N_v\otimes\det V_v)^{-1}\Bigr)
\notag\\
&\label{uy5}\hspace{40pt}
\otimes\bigotimes_{v>0,\,v\equiv0\,{\rm mod}\,(n_j)}\Lambda V_v\ ,
\\
S(U_2,L_2)'&=S\Bigl(TX_H,L_2\otimes\bigotimes_{v>0,v\equiv0\ {\rm mod}(n_j)}
(\det N_v)^{-1}\Bigr.
\notag\\
&\label{uy6}\hspace{20pt}
\Bigl.\otimes\bigotimes_{v>0,\,v\equiv \tfrac{n_j}{2}\,{\rm mod}\,(n_j)}
(\det V_v)^{-1}\Bigr)
\otimes\bigotimes_{v>0,\,v\equiv\tfrac{n_j}{2}\,{\rm mod}\,(n_j)}\Lambda V_v\ .
\end{align}
Then by (\ref{liu17}) and (\ref{liu18}), for $i=1, 2$, we have the following isomorphisms of
Clifford modules over $X_H$,
\begin{equation}\label{eg58}
S(U_i,L_i)\simeq S(U_i,L_i)'\otimes\bigotimes_{v>0,\,v\equiv0\,{\rm mod}\,(n_j)}\Lambda N_v.
\end{equation}
We define the $\mathbb{Z}_2$-gradings on $S(U_i,L_i)'$ ($i=1,2$) induced by the $\mathbb{Z}_2$-gradings on
$S(U_i,L_i)$ ($i=1,2$) and on $\otimes_{v>0,\,v\equiv0\,{\rm mod}\,(n_j)}\Lambda N_v$ such that the isomorphisms
(\ref{eg58}) preserve the $\mathbb{Z}_2$-gradings.

As in \cite[pp. 952]{MR1870666}, we introduce formally the following $\mathbb{Z}/k$ complex line bundles over $X_H$,
\begin{equation*}\label{eg59}
\begin{split}
L_1'&=\Bigl(L_1^{-1}\otimes\negthickspace\bigotimes_{\negthickspace v>0,\,v\equiv0\,{\rm mod}\,(n_j)}
(\det N_v\otimes\det V_v)
\bigotimes_{v>0}\,(\det N_v\otimes\det V_v)^{-1}\otimes K_X\Bigr)^{\frac{1}{2}},
\\
L_2'&=\Bigl(\!L_2^{-1}\otimes\negthickspace\bigotimes_{\negthickspace v>0,\,v\equiv0\,{\rm mod}\,(n_j)}\negthickspace\det N_v
\bigotimes_{\negthickspace v>0,\,v\equiv\frac{n_j}{2}\,{\rm mod}\,(n_j)}\negthickspace\det V_v
\bigotimes_{v>0}\,(\det N_v)^{-1}\otimes K_X\Bigr)^{\frac{1}{2}}.
\end{split}
\end{equation*}
In fact, from (\ref{liu17}), (\ref{liu18}), Lemma \ref{eg54} and the assumption that $V$ is spin,
one verifies easily that $c_1(L_i'^2)=0 \mod (2)$ for $i=1,2$, which implies that $L_1'$ and $L_2'$ are
well defined $\mathbb{Z}/k$ complex line bundles over $X_H$ (cf. Section \ref{sec2.1}).

Then by \eqref{uy5}, \eqref{uy6} and the definitions of $L_1$, $L_2$, $L'_1$ and $L'_2$, we get the following identifications of $\mathbb{Z}/k$ Clifford modules over $X_H$ (cf. \cite[(4.19)]{MR1870666}),
\begin{align}\label{eg60}
S(U_1,L_1)'\otimes L_1'&=S\bigl(TX_H,K_X\otimes_{v>0}(\det N_v)^{-1}\bigr)\otimes S\bigl(V_0^{\mathbb{R}},
\otimes_{v>0}(\det V_v)^{-1}\bigr)
\notag\\
&\qquad\otimes\bigotimes_{v>0,\,v\equiv0\,{\rm mod}\,(n_j)}\Lambda(V_v)\ ,
\\
\label{uy30}
S(U_2,L_2)'\otimes L_2'&=S\bigl(TX_H,K_X\otimes_{v>0}(\det N_v)^{-1}\bigr)\otimes\bigotimes_{v>0,\,v\equiv \frac{n_j}{2}
\,{\rm mod}\,(n_j)}\Lambda(V_v)\ .
\end{align}

\begin{lem}\label{eg61}
{\rm (cf. \cite[Lemma 4.3]{MR1870666})} Let us write
\begin{equation*}\label{eg62}
\begin{split}
L(\beta_j)_1&=L_1'\otimes\bigotimes_{v>0}(\det N_v)^{\bigl[\tfrac{p_j}{n_j}v\bigr]+(p-1)v+1}
\otimes\bigotimes_{v>0}(\det \overline{V}_v)^{\bigl[\tfrac{p_j}{n_j}v\bigr]+(p-1)v}
\\
&\qquad\otimes\bigotimes_{v>0,\,v\equiv0\,{\rm mod}(n_j)}
(\det N_v)^{-1}\otimes\bigotimes_{v<0}(\det W_v)^{\bigl[-\tfrac{p_j}{n_j}v\bigr]-(p-1)v}
\\
&\qquad\otimes\bigotimes_{v>0}(\det \overline{W}_v)^{\bigl[\tfrac{p_j}{n_j}v\bigr]+(p-1)v+1}
\otimes\bigotimes_{v>0,\,v\equiv0\,{\rm mod}(n_j)}
(\det \overline{W}_v)^{-1}\ ,
\\
L(\beta_j)_2&=L_2'\otimes\bigotimes_{v>0}(\det N_v)^{\bigl[\tfrac{p_j}{n_j}v\bigr]+(p-1)v+1}
\otimes\bigotimes_{v>0}(\det \overline{V}_v)^{\bigl[\tfrac{p_j}{n_j}v+\tfrac{1}{2}\bigr]+(p-1)v}
\\&\qquad
\otimes\bigotimes_{v>0,\,v\equiv0\,{\rm mod}(n_j)}
(\det N_v)^{-1}\otimes\bigotimes_{v<0}(\det W_v)^{\bigl[-\tfrac{p_j}{n_j}v\bigr]-(p-1)v}
\\&\qquad
\otimes\bigotimes_{v>0}(\det \overline{W}_v)^{\bigl[\tfrac{p_j}{n_j}v\bigr]+(p-1)v+1}\otimes
\bigotimes_{v>0,\,v\equiv0\,{\rm mod}(n_j)}(\det \overline{W}_v)^{-1}\ .
\end{split}
\end{equation*}
Then $L(\beta_j)_1$ and $L(\beta_j)_2$ can be extended naturally to $\mathbb{Z}/k$ $G_y\times S^1$-equivariant complex
line bundles over $X(n_j)$ which we will still denote by $L(\beta_j)_1$ and $L(\beta_j)_2$  respectively.
\end{lem}

Now we compare the $\mathbb{Z}_2$-gradings in (\ref{eg60}). Set
\begin{equation}\label{uy7}
\begin{split}
\Delta(n_j,N)&=\sum_{\tfrac{n_j}{2}<v'<n_j}\sum_{0<v,\,v\equiv v'\,{\rm mod}\,(n_j)}
\dim N_v+o\Bigl(N(n_j)_{\tfrac{n_j}{2}}^{\mathbb{R}}\Bigr)\ ,
\\
\Delta(n_j,V)&=\sum_{\tfrac{n_j}{2}<v'<n_j}\sum_{0<v,\,v\equiv v'\,{\rm mod}\,(n_j)}
\dim V_v+o\Bigl(V(n_j)_{\tfrac{n_j}{2}}^{\mathbb{R}}\Bigr)\ ,
\end{split}
\end{equation}
with $o(N(n_j)_{\tfrac{n_j}{2}}^{\mathbb{R}})$ (resp. $o(V(n_j)_{\tfrac{n_j}{2}}^{\mathbb{R}})$ equals 0 or 1,
depending on whether the given orientation on $N(n_j)_{\tfrac{n_j}{2}}^{\mathbb{R}}$
(resp. $V(n_j)_{\tfrac{n_j}{2}}^{\mathbb{R}}$) agrees or disagrees with the complex
orientation of $\oplus_{v>0,\,v\equiv \tfrac{n_j}{2}\,{\rm mod}\,(n_j)} N_v$
(resp. $\oplus_{v>0,\,v\equiv \tfrac{n_j}{2}\,{\rm mod}\,(n_j)} V_v$).

By \cite[pp. 953]{MR1870666}, we know that for the $\mathbb{Z}_2$-gradings induced by $\tau_s$,
the differences of the $\mathbb{Z}_2$-gradings of (\ref{eg60}) and \eqref{uy30} are both $(-1)^{\Delta(n_j,N)}$;
for the $\mathbb{Z}_2$-gradings induced by $\tau_e$,
the difference of the $\mathbb{Z}_2$-gradings of (\ref{eg60}) (resp. \eqref{uy30}) is
$(-1)^{\Delta(n_j,N)+\Delta(n_j,V)}$ (resp. $\,(-1)^{\Delta(n_j,N)+o\bigl(V(n_j)_{n_j/2}^{\mathbb{R}}\bigr)}\,$).

To simplify the notations, we introduce the same functions as in \cite[(4.30)]{MR1870666}, which are locally constant on $X_H$,
\begin{equation}
\begin{split}
\varepsilon(W)&=-\frac{1}{2}\sum_{v>0}(\dim W_v)\cdot\Bigl(\bigl(\bigl[\tfrac{p_j}{n_j}v\bigr]+(p-1)v\bigr)
\bigl(\bigl[\tfrac{p_j}{n_j}v\bigr]+(p-1)v+1\bigr)\Bigr.
\\
&\hspace{60pt}
-\Bigl.(\tfrac{p_j}{n_j}v+(p-1)v)\bigl(2\bigl(\bigl[\tfrac{p_j}{n_j}v\bigr]+(p-1)v\bigr)+1\bigr)\Bigr)
\\
&-\frac{1}{2}\sum_{v<0}(\dim W_v)\cdot\Bigl(\bigl(\bigl[-\tfrac{p_j}{n_j}v\bigr]-(p-1)v\bigr)
\bigl(\bigl[-\tfrac{p_j}{n_j}v\bigr]-(p-1)v+1\bigr)\Bigr.
\\
&\hspace{60pt}
+\Bigl.\bigl(\tfrac{p_j}{n_j}v+(p-1)v\bigr)\bigl(2\bigl(\bigl[-\tfrac{p_j}{n_j}v\bigr]-(p-1)v\bigr)+1\bigr)\Bigr)\ ,
\end{split}
\end{equation}
\begin{align}
\varepsilon_1&=\frac{1}{2}\sum_{v>0}(\dim N_v-\dim V_v)\Big(\bigl(\bigl[\tfrac{p_j}{n_j}v\bigr]+(p-1)v\bigr)
\bigl(\bigl[\tfrac{p_j}{n_j}v\bigr]+(p-1)v+1\bigr)\Bigr.
\notag\\
&\hspace{60pt}
\Bigl.-\bigl(\tfrac{p_j}{n_j}v+(p-1)v\bigr)\bigl(2\bigl(\bigl[\tfrac{p_j}{n_j}v\bigr]+(p-1)v\bigr)+1\bigr)\Bigr)\ ,
\end{align}
\begin{equation}
\begin{split}
\varepsilon_2&=\frac{1}{2}\sum_{v>0}(\dim N_v)\cdot\Bigl(\bigl(\bigl[\tfrac{p_j}{n_j}v\bigr]+(p-1)v\bigr)
\bigl(\bigl[\tfrac{p_j}{n_j}v\bigr]+(p-1)v+1\bigr)\Bigr.
\\
&\hspace{70pt}
\Bigl.-(\tfrac{p_j}{n_j}v+(p-1)v)\bigl(2\bigl(\bigl[\tfrac{p_j}{n_j}v\bigr]+(p-1)v\bigr)+1\bigr)\Bigr)
\\
&\quad-\frac{1}{2}\sum_{v>0}(\dim V_v)\cdot\Bigl(\bigl(\bigl[\tfrac{p_j}{n_j}+\tfrac{1}{2}\bigr]+(p-1)v\bigl)^2\Bigr.
\\
&\hspace{70pt}
\Bigl.-2\bigl(\tfrac{p_j}{n_j}v+(p-1)v\bigr)\bigl(\bigl[\tfrac{p_j}{n_j}+\tfrac{1}{2}\bigr]+(p-1)v\bigr)\Bigr)\ .
\end{split}
\end{equation}
As in \cite[(2.23)]{MR1870666}, for $0\leq j\leq J_0$, we set
\begin{align}\label{eg64}
e(p,\beta_j,N)&=\frac{1}{2}\sum_{v>0}(\dim N_v)\cdot\bigl(\bigl[\tfrac{p_j}{n_j}v\bigr]+(p-1)v\bigr)
\bigl(\bigl[\tfrac{p_j}{n_j}v\bigl]+(p-1)v+1\bigl)\ ,
\notag\\
d'(p,\beta_j,N)&=\sum_{v>0}(\dim N_v)\cdot\bigl(\bigl[\tfrac{p_j}{n_j}v\bigr]+(p-1)v\bigr)\ .
\end{align}
Then $e(p,\beta_j,N)$ and $d'(p,\beta_j,N)$ are locally constant functions on $X_H$. In particular, we have
\begin{equation}\label{eg65}
\begin{split}
e(p,\beta_0,N)&=\frac{1}{2}(p-1)^2e(N)+\frac{1}{2}(p-1)d'(N)\ ,
\\
e(p,\beta_{J_0},N)&=\frac{1}{2}p^2e(N)+\frac{1}{2}p\,d'(N)\ ,
\\
d'(p,\beta_{J_0},N)&=d'(p+1,\beta_{0},N)=p\,d'(N)\ .
\end{split}
\end{equation}

By Proposition \ref{eg44}, (\ref{eg60}) and Lemma \ref{eg62}, we deduce an analogue of Proposition \ref{eg25}.
\begin{prop}\label{eg66}{\rm (cf. \cite[Proposition 4.2]{MR1870666})}
For $i=1,2$, the  $\mathbb{Z}/k$ $G_y$-equivariant isomorphisms of complex vector bundles over $X_H$,
\begin{equation*}\label{eg67}
\begin{split}
&r_{i1}:\,S(TX_H,K_X\otimes_{v>0}(\det N_v)^{-1})\otimes(K_W\otimes K^{-1}_X)^{1/2}
\\&\hspace{60pt}
\otimes \mathcal{F}_{p,j-1}(X)
\otimes F_V^i\otimes Q^1(W)
\\
&\qquad\longrightarrow S(U_i,L_i)\otimes(K_W\otimes K^{-1}_X)^{1/2}\otimes \mathcal{F}(\beta_j)\otimes F_V^i(\beta_j)
\\&\hspace{60pt}
\otimes Q_W(\beta_j)\otimes L(\beta_j)_i\otimes\bigotimes_{v>0,\,v\equiv0\,{\rm mod}\,(n_j)}{\rm Sym}\,(\overline{N}_{v,0})\ ,
\\
&r_{i2}:\,S(TX_H,K_X\otimes_{v>0}(\det N_v)^{-1})\otimes(K_W\otimes K^{-1}_X)^{1/2}
\\&\hspace{60pt}
\otimes \mathcal{F}_{p,j}(X)
\otimes F_V^i\otimes Q^1(W)\\
&\qquad\longrightarrow S(U_i,L_i)\otimes(K_W\otimes K^{-1}_X)^{1/2}\otimes \mathcal{F}(\beta_j)\otimes F_V^i(\beta_j)
\\&\hspace{60pt}
\otimes Q_W(\beta_j)\otimes L(\beta_j)_i\otimes
\bigotimes_{v>0,\,v\equiv0\,{\rm mod}\,(n_j)}
\bigl({\rm Sym}\,(N_{v,0})\otimes \det N_v\bigr)
\end{split}
\end{equation*}
have the following properties:
\begin{enumerate}[{\rm (i)}]
\item for $i=1,2$, $\gamma=1,2$,
\begin{equation}\label{eg68}
\begin{split}
r_{i\gamma}^{-1}\cdot \textbf{J}_H\cdot r_{i\gamma}&=\textbf{J}_H\ ,
\\
r_{i\gamma}^{-1}\cdot P\cdot r_{i\gamma}&=P+\Bigl(\,\frac{p_j}{n_j}+(p-1)\Bigr)\textbf{J}_H+\varepsilon_{i\gamma}\ ,
\end{split}
\end{equation}
where $\varepsilon_{i1}=\varepsilon_i+\varepsilon(W)-e(p,\beta_{j-1},N)$,\ $\varepsilon_{i2}=\varepsilon_i+\varepsilon(W)-e(p,\beta_{j},N)$.

\item Recall that $o\bigl(V(n_j)_{\tfrac{n_j}{2}}^{\mathbb{R}}\bigr)$ is defined in \eqref{uy7}. Let
\begin{equation*}
\begin{split}
\mu_1&=-\sum_{v>0}\bigl[\tfrac{p_j}{n_j}v\bigr]\dim V_v+\Delta(n_j,N)+\Delta(n_j,V) \mod (2),
\\
\mu_2&=-\sum_{v>0}\bigl[\tfrac{p_j}{n_j}v+\tfrac{1}{2}\bigr]\dim V_v+\Delta(n_j,N)
+o\bigl(V(n_j)_{\tfrac{n_j}{2}}^{\mathbb{R}}\bigr) \mod (2),
\\
\mu_3&=\Delta(n_j,N) \mod (2),
\\
\mu_4&=\sum_{v}\bigl(\bigl[\tfrac{p_j}{n_j}v\bigr]+(p-1)v\bigr)\dim W_v+\dim W+\dim W(n_j)_0 \mod (2).
\end{split}
\end{equation*}
Then for $i=1,2$, $\gamma=1,2$, we have
\begin{equation}\label{eg71}
\begin{split}
r_{i\gamma}^{-1}\tau_{e}r_{i\gamma}&=(-1)^{\mu_i}\tau_{e}\ ,\quad r_{i\gamma}^{-1}\tau_{s}r_{i\gamma}
=(-1)^{\mu_3}\tau_{s}\ ,
\\
r_{i\gamma}^{-1}\tau_{1}r_{i\gamma}&=(-1)^{\mu_4}\tau_{1}\ .
\end{split}
\end{equation}
\end{enumerate}
\end{prop}

\subsection{A proof of Theorem \ref{eg31}}\label{sec5.3}
Let $X'$ be a connected component  of $X(n_j)$. By \cite[Lemmas 4.4, 4.5, 4.6]{MR1870666}, we know that for
$i=1,2$, $k=1,2,3$, the following functions are independent on the connected components of $X_H$ in $X'$,
\begin{equation*}
\begin{split}
&\varepsilon_i+\varepsilon(W) \mod(2)\ ,\quad d'(p,\beta_j,N)+\mu_k+\mu_4 \mod(2)\ ,
\\
&\sum_{v>0}\bigl[\tfrac{p_j}{n_j}v\bigr]\dim V_v+\Delta(n_j,V) \mod(2)\ ,
\\
&\sum_{v>0}\bigl[\tfrac{p_j}{n_j}v+\tfrac{1}{2}\bigr]\dim V_v+o\bigl(V(n_j)_{\tfrac{n_j}{2}}^{\mathbb{R}}\bigr) \mod(2)\ ,
\end{split}
\end{equation*}
which implies that $d'(p,\beta_{j-1},N)+\sum_{0<v}\dim N_v+\mu_k+\mu_4 \mod(2)$ ($k=1,2,3$) are constant functions on each connected
components of $X(n_j)$.

By \eqref{uy8}, \eqref{uy9}, \eqref{uy10}, (\ref{eg52}) and Lemma \ref{eg61}, we know that the Dirac operator
$D^{X(n_j)}\otimes\mathcal{F}(\beta_j)\otimes F_V^i(\beta_j)
\otimes Q_W(\beta_j)\otimes L(\beta_j)_i$ ($i=1,2$) is well defined on $X(n_j)$.
Observe that the two equalities in Theorem \ref{liu28} are both compatible with the $G_y$ action.
Thus, by using Proposition \ref{eg66} and applying both the first and the second
equalities of Theorem \ref{liu28} to each connected component of
$X(n_j)$ separately, we deduce that for $i=1,2$, $1\leq j\leq J_0$,
$m\in \tfrac{1}{2}\mathbb{Z}$, $1\leq \ell<N$,
$h\in \mathbb{Z}$, $\tau=\tau_{e1}$ or $\tau_{s1}$,
\begin{align}\label{uy11}
&\sum_{\alpha}(-1)^{d'(p,\beta_{j-1},N)+\sum_{v>0}\dim N_v}
\APS_{\tau}\Bigl(D^{X_{H,\alpha}}\otimes(K_W\otimes K^{-1}_X)^{1/2}\Bigr.
\notag\\
&\hspace{40pt}
\Bigl.\otimes\mathcal{F}_{p,j-1}(X)\otimes F_V^i\otimes\ Q^1(W),m+e(p,\beta_{j-1},N),\ell,h\Bigr)
\notag\\
&\equiv\sum_{\beta}(-1)^{d'(p,\beta_{j-1},N)+\sum_{v>0}\dim N_v+\mu}
\APS_{\tau}\Bigl(D^{X(n_j)}\otimes(K_W\otimes K^{-1}_X)^{1/2}\Bigr.
\notag\\
&\hspace{60pt}
\otimes\mathcal{F}(\beta_j)\otimes F_V^i(\beta_j)
\otimes Q_W(\beta_j)\otimes L(\beta_j)_i,m+\varepsilon_i
\notag\\
&\hspace{120pt}
\Bigl.+\varepsilon(W)+(\tfrac{p_j}{n_j}+(p-1))h,\ell,h\Bigr)
\notag\\
&\equiv\sum_{\alpha}(-1)^{d'(p,\beta_{j},N)+\sum_{v>0}\dim N_v}
\APS_{\tau}\Bigl(D^{X_{H,\alpha}}\otimes(K_W\otimes K^{-1}_X)^{1/2}\Bigr.
\notag\\
&\hspace{40pt}
\Bigl.\otimes\mathcal{F}_{p,j}(X)\otimes F_V^i\otimes Q^1(W),m+e(p,\beta_{j},N),\ell,h\Bigr) \mod k\mathbb{Z}\ ,
\end{align}
where $\sum_{\beta}$ means the sum over all the connected components of $X(n_j)$. In \eqref{uy11}, if $\tau=\tau_{s1}$, then
$\mu=\mu_3+\mu_4$; if $\tau=\tau_{e1}$, then $\mu=\mu_i+\mu_4$.
Combining \eqref{eg65} with \eqref{uy11}, we get \eqref{eg32}.

The proof of Theorem \ref{eg31} is completed.

\vspace{3mm}\textbf{Acknowledgements}\ \ The authors wish to thank
Professors Daniel S. Freed, Xiaonan Ma and Weiping Zhang for their helpful
discussions.


\end{document}